\documentclass[11pt]{article}
\usepackage{amssymb,amsmath,enumerate}
\newcommand{\ff}{\mathbf{f}}
\newcommand{\uu}{\mathbf{u}}
\newcommand{\vv}{\mathbf{v}}
\newcommand{\ww}{\mathbf{w}}
\newcommand{\xx}{\mathbf{x}}
\newcommand{\yy}{\mathbf{y}}
\newcommand{\ee}{\mathbf{e}}
\newcommand{\GG}{\mathbf{G}}
\newcommand{\oomega}{{\boldsymbol\omega}}

\newcommand{\cA}{\mathcal{A}}

\newcommand{\cG}{\mathcal{G}}
\newcommand{\cI}{\mathcal{I}}
\newcommand{\cL}{\mathcal{L}}
\newcommand{\cM}{\mathcal{M}}
\newcommand{\cO}{\mathcal{O}}

\newcommand{\C}{\mathbb{C}}
\newcommand{\N}{\mathbb{N}}
\newcommand{\PP}{\mathbb{P}}
\newcommand{\R}{\mathbb{R}}

\newcommand{\X}{\mathbb{X}}
\newcommand{\Z}{\mathbb{Z}}
\newcommand{\dd}{\,{\rm d}}
\newcommand{\D}{{\rm d}}

\renewcommand{\div}{\mathop{\mathrm{div}}}
\newcommand{\curl}{\mathop{\mathrm{curl}}}

\renewcommand{\:}{\thinspace :}

\newcommand{\weakto}{\rightharpoonup}

\renewcommand{\Re}{\mathop{\mathrm{Re}}}
\renewcommand{\Im}{\mathop{\mathrm{Im}}}

\numberwithin{equation}{section}
\newtheorem{thm}{Theorem}[section]

\newtheorem{prop}[thm]{Proposition}
\newtheorem{lem}[thm]{Lemma}
\newtheorem{rem}[thm]{Remark}

\newcommand{\QED}{\mbox{}\hfill$\Box$}

\begin{document}

\title{Existence and stability of viscous vortices%
\footnote{To appear in ``Handbook of Mathematical Analysis 
in Mechanics of Viscous Fluids'', Y. Giga and A. Novotn\'y Ed., Springer. 
The final publication is available at http://www.springerlink.com.}}

\author{
\null\\
{\bf Thierry Gallay}\\ 
Universit\'e Grenoble Alpes\\
Institut Fourier\\
UMR CNRS 5582\\
100 rue des Maths\\
38610 Gi\`eres, France\\
{\tt Thierry.Gallay@univ-grenoble-alpes.fr}
\and
\\
{\bf Yasunori Maekawa}\\    
Department of Mathematics\\
Graduate School of Science\\ 
Kyoto University\\ 
Kitashirakawa Oiwake-cho, Sakyo-ku\\ 
Kyoto 606-8502, Japan\\
{\tt maekawa@math.kyoto-u.ac.jp}}

\date{}

\maketitle

\begin{abstract}
Vorticity plays a prominent role in the dynamics of incompressible
viscous flows. In two-dimensional freely decaying turbulence, after a
short transient period, evolution is essentially driven by
interactions of viscous vortices, the archetype of which is the
self-similar Lamb-Oseen vortex. In three dimensions, amplification of
vorticity due to stretching can counterbalance viscous dissipation and
produce stable tubular vortices. This phenomenon is illustrated in a
famous model originally proposed by Burgers, where a straight vortex
tube is produced by a linear uniaxial strain field. In real flows
vortex lines are usually not straight, and can even form closed curves,
as in the case of axisymmetric vortex rings which are very common in
nature and in laboratory experiments. The aim of this chapter is 
to review a few rigorous results concerning existence and stability 
of viscous vortices in simple geometries.  
\end{abstract}

\section{Introduction}\label{sec0}

Since the pioneering work of Helmholtz \cite{vHe}, vorticity has been
widely recognized as a quantity of fundamental importance in fluid
dynamics, especially for turbulent flows. According to a famous quote
by K\"uchermann \cite{Ku}, ``vortices are the sinews and muscles of
fluid motions''. Intuitively, vorticity describes the local rotation
of fluid particles at a given point. In the Eulerian representation,
if $\uu(\xx,t)$ denotes the velocity of the fluid at point $\xx =
(x_1,x_2,x_3) \in \R^3$ and time $t \in \R$, the vorticity is the vector
$\oomega(\xx,t) = \curl \uu(\xx,t) = \nabla \wedge \uu(\xx,t)$. 
Under the evolution given by the Navier-Stokes equations, the vorticity 
satisfies
\begin{equation}\label{3DNSom}
  \partial_t \oomega(\xx,t) + (\uu(\xx,t),\nabla)\oomega(\xx,t) 
  - (\oomega(\xx,t),\nabla)\uu(\xx,t) \,=\, \nu\Delta \oomega(\xx,t)\,, 
\end{equation}
where $\nu > 0$ is the kinematic viscosity of the fluid, i.e. the
ratio of the viscosity to the fluid density. In the incompressible
case considered here, the velocity field satisfies $\div \uu(\xx,t) =
0$ and is thus entirely determined by the vorticity distribution up to
an irrotational flow. The Biot-Savart law is a reconstruction formula
that expresses $\uu$ in terms of $\oomega$, depending on the geometry
of the fluid domain and the boundary conditions. In the whole space
$\R^3$, if the vorticity distribution is sufficiently localized, the
Biot-Savart formula reads
\begin{equation}\label{BS3D}
  \uu(\xx,t) \,=\, -\frac{1}{4\pi}\int_{\R^3}\frac{(\xx-\yy)\wedge 
  \oomega(\yy,t)}{|\xx-\yy|^3} \dd \yy\,.
\end{equation}

In viscous fluids, vorticity is usually created within boundary layers
near walls or interfaces, or in the vicinity of a stirring
device. Once produced, vorticity can be substantially amplified by the
local strain in the fluid, through a genuinely three-dimensional
mechanism that is often referred to as ``vortex stretching''. A Taylor
expansion of the velocity field at a given point $\xx_0$ reveals that
\[
  \uu(\xx,t) \,=\, \uu(\xx_0,t) + \frac12 \oomega(\xx_0,t) \wedge
  (\xx-\xx_0) + (D\uu(\xx_0,t))(\xx-\xx_0) + \cO(|\xx-\xx_0|^2)\,,
\]
where $D\uu = \frac12 \bigl((\nabla\uu) + (\nabla\uu)^\top\bigr)$ is
the deformation tensor, whose eigenvalues $\gamma_1, \gamma_2,
\gamma_3$ are called the principal strains at $\xx_0$.
Incompressibility implies that $\gamma_1 + \gamma_2 + \gamma_3 = 0$,
so that two generic situations may occur. If two principal strains
(say, $\gamma_1$ and $\gamma_2$) are negative and the third one is
positive, vorticity gets amplified at $\xx_0$ in the direction of the
principal strain axis corresponding to $\gamma_3$. That stretching
mechanism can compensate the viscous dissipation and result in the
formation of stable vortex filaments, a typical example being the
Burgers vortex \cite{Bu} which will be studied in Section~\ref{sec3}.
In contrast, if two principal strains are positive at $\xx_0$, the
stretching effect leads to the formation of vortex sheets, which are
also commonly observed in turbulent flows although they undergo the
Kelvin-Helmholtz instability at high Reynolds numbers. Vortex sheets
play a prominent role in interfacial motion and boundary layer theory,
and the interested reader is referred to the chapter entitled
``The Inviscid Limit and Boundary Layers for the Navier-Stokes Flows'' 
for further information.

In the present chapter, emphasis is put on vortex tubes or filaments,
for which vorticity is essentially concentrated along a curve with no
endpoints in the fluid. In general, the curve will evolve with time,
because it is advected by the flow. According to Helmholtz's first
law, the \emph{total circulation} of such a vortex filament is
constant along its length, and is also independent of time as long as
the viscous effects can be neglected. This very important quantity,
often denoted by $\Gamma$, can be defined as the flux of the vorticity
vector through any cross section of the vortex tube, or equivalently
(in view of Stokes' theorem) as the circulation of the velocity along
any closed curve enclosing that tube. The ratio $\alpha = \Gamma/\nu$
of the total circulation to the kinematic viscosity is a dimensionless
quantity, sometimes referred to as the circulation Reynolds number,
which measures the strength of the vortex and plays a crucial role in
stability issues.

Since vortex filaments have no endpoints, they must either extend to
the fluid boundary, or to infinity, or form closed curves.  In the
simple situation, already considered in \cite{vHe}, where all vortex
lines are straight and parallel to each other, the velocity and
vorticity fields take the particular form
\begin{equation}\label{velvort2D}
  \uu(\xx,t) \,=\, \begin{pmatrix} u_1(x_1,x_2,t) \\ u_2(x_1,x_2,t)
  \\ 0 \end{pmatrix}\,, \qquad 
  \oomega(\xx,t) \,=\, \begin{pmatrix} 0 \\ 0 \\ \omega(x_1,x_2,t)
  \end{pmatrix}\,,
\end{equation}
where $x = (x_1,x_2) \in \R^2$ and $\omega = \partial_1 u_2 - 
\partial_2 u_1$. Here the coordinates have been chosen so that the
third axis coincides with the direction of the vortex filaments.
The evolution equation for the scalar vorticity $\omega(x,t)$ is
\begin{equation}\label{2DNSom}
  \partial_t \omega(x,t) + u(x,t)\cdot\nabla\omega(x,t) 
  \,=\, \nu\Delta \omega(x,t)\,, 
\end{equation}
where $u = (u_1,u_2)$ satisfies $\partial_1 u_1 + \partial_2 u_2 = 0$. 
If the vorticity distribution is sufficiently localized, the 
two-dimensional velocity field $u(x,t)$ is given by the 
2D Biot-Savart law
\begin{equation}\label{BS2D}
  u(x,t) \,=\, \frac{1}{2\pi}\int_{\R^2}\frac{(x-y)^\perp}{|x-y|^2} 
  \,\omega(y,t)\dd y\,,
\end{equation}
where $x^\perp = (-x_2,x_1)$ and $|x|^2 = x_1^2 + x_2^2$. 
Eq.~\eqref{2DNSom} is just an advection-diffusion equation for the
scalar quantity $\omega$, hence (by the maximum principle) no
amplification of vorticity can occur in the two-dimensional case. As a
consequence, all localized vortex structures will eventually spread
out and decay, since there is nothing to counterbalance the effect of
viscosity. A typical example is provided by the \emph{Lamb-Oseen vortex}, 
an exact self-similar solution to \eqref{2DNSom} of the form
\begin{equation}\label{Osdef}
  \omega(x,t) \,=\, \frac{\Gamma}{\nu t}\,G\Bigl(\frac{x}{\sqrt{\nu t}}
  \Bigr)\,,\qquad u(x,t) \,=\, \frac{\Gamma}{\sqrt{\nu t}}\,
  v^G\Bigl(\frac{x}{\sqrt{\nu t}}\Bigr)\,,
\end{equation}
where the vorticity and velocity profiles are explicitly given by
\begin{equation}\label{Gdef}
  G(\xi) \,=\, \frac{1}{4\pi}\,e^{-|\xi|^2/4}\,, \qquad 
  v^G(\xi) \,=\, \frac{1}{2\pi}\,\frac{\xi^\perp}{|\xi|^2}
  \Bigl(1 - e^{-|\xi|^2/4}\Bigr)\,, \qquad \xi \in \R^2\,.
\end{equation}
Note that $\int_{\R^2} G(\xi)\dd\xi = 1$, so that $\int_{\R^2}
\omega(x,t)\dd x = \Gamma$ for all $t > 0$, in agreement with the
general definition of the total circulation $\Gamma$. The Lamb-Oseen
vortex plays a distinguished role in the dynamics of the
two-dimensional vorticity equation \eqref{2DNSom}, for two main
reasons. First, it deserves the name of fundamental solution, in the
sense that it is the unique solution of \eqref{2DNSom} with initial
data $\omega_0 = \Gamma\delta_0$, where $\delta_0$ denotes the Dirac
measure at the origin. Next, it describes to leading order the
long-time asymptotics of all solutions of \eqref{2DNSom} with
integrable initial data and nonzero circulation \cite{GW1}. If
self-similar variables are used, the Lamb-Oseen vortex becomes a
stationary solution of some rescaled equation, and its stability
properties can then be studied using spectral theory and other
standard techniques. This analysis is presented in Section~\ref{sec1}
below, and serves as a model for further existence and stability
results in more complex situations.

Another relatively simple and mathematically tractable situation is
the axisymmetric case without swirl, where the velocity field is
invariant under rotations about a given axis, and under reflections by
any plane containing the axis. Here all vortex lines are circles
centered on the symmetry axis and normal to it. Using cylindrical
coordinates $(r,\theta,z)$, so that $r$ represents the distance to the
symmetry axis and $z$ the position along the axis, the velocity and
vorticity fields are given by
\begin{equation}\label{velvortaxi}
  \uu(\xx,t) \,=\, u_r(r,z,t)\,\ee_r + u_z(r,z,t)\,\ee_z\,,
  \qquad \oomega(\xx,t) \,=\, \omega_\theta(r,z,t)\,\ee_\theta\,,
\end{equation}
where $\ee_r, \ee_\theta, \ee_z$ denote unit vectors in the radial, 
toroidal, and vertical directions, respectively. As in the
two-dimensional case, the vorticity vector has only one nonzero
component $\omega_\theta$, which satisfies the evolution equation
\begin{equation}\label{Axiom}
  \partial_t \omega_\theta + u\cdot\nabla\omega_\theta -
  \frac{u_r}{r}\omega_\theta \,=\, \nu\Bigl(\Delta \omega_\theta
  - \frac{\omega_\theta}{r^2}\Bigr)\,,
\end{equation}
where $u\cdot\nabla = u_r\partial_r + u_z\partial_z$ and $\Delta =
\partial_r^2 + \frac{1}{r}\partial_r + \partial_z^2$ denotes the
Laplace operator in cylindrical coordinates. The velocity 
$u = (u_r,u_z)$ can be expressed in terms of the axisymmetric
vorticity $\omega_\theta$ by solving the linear elliptic system
\begin{equation}\label{AxiBS}
  \partial_r u_r + \frac{1}{r} u_r + \partial_z u_z \,=\, 0\,,
  \qquad \partial_z u_r - \partial_r u_z \,=\, \omega_\theta\,,
\end{equation}
in the half-plane $\Omega = \{(r,z) \in \R^2\,|\, r > 0\,,~z \in
\R\}$, with boundary conditions $u_r = \partial_r u_z = 0$ at $r = 0$.
Explicit formulas for the axisymmetric Biot-Savart law exist, see
e.g. \cite{FS,GS1}, but are more involved than in the two-dimensional
case. The analogue of the Lamb-Oseen vortex for axisymmetric flows is
the solution of \eqref{Axiom} with a vortex filament as initial data.
This means that $\omega_\theta(\cdot,\cdot,0) = \Gamma\delta_{(\bar r,
\bar z)}$ where $\delta_{(\bar r,\bar z)}$ denotes the Dirac measure
located at some point $(\bar r,\bar z) \in \Omega$. Existence of a
global solution to \eqref{Axiom} with such initial data was recently
shown by Feng and \v{S}ver\'ak \cite{FS}, and uniqueness can be
established using, in particular, the approach presented in
Section~\ref{sec1}, see \cite{GS2}. The reader is referred to
Section~\ref{sec2} below for up-to-date results on existence of
axisymmetric vortex rings.

The third and final case considered here is a famous model for vortex
filaments in turbulent flows, originally proposed by Burgers
\cite{Bu}. It is assumed that the velocity field has the form
$\uu(\xx,t) = \uu_s(\xx) + \vv(\xx,t)$, where $\uu_s(\xx)$ is a 
stationary straining flow of the form
\begin{equation}\label{defstrain}
  \uu_s(\xx) \,=\, \begin{pmatrix} \gamma_1 x_1 \\ \gamma_2 x_2  \\ 
  \gamma_3 x_3\end{pmatrix} \,=\, \cM \xx\,, \qquad
   \cM \,=\, \begin{pmatrix} \gamma_1 & 0 & 0\\ 0 & \gamma_2 & 0 \\ 
  0 & 0 & \gamma_3\end{pmatrix}\,,
\end{equation}
where $\gamma_1 + \gamma_2 + \gamma_3 = 0$ and $\gamma_1, \gamma_2 <
0$, $\gamma_3 > 0$. According to the discussion above, the strain
\eqref{defstrain} describes to leading order the deformation rate of
any smooth, incompressible velocity field near the origin, at a given
time. Burgers' model is crude in the sense that it assumes that the
strain $\uu_s(\xx)$ is independent of time and extends all the way to
infinity in space, which is certainly not realistic in turbulent
flows. Nevertheless, the model is interesting because it clearly
illustrates the vortex stretching effect, which in the present case
produces a family of stationary solutions that can be compared 
with observations in experiments.

If $\uu(\xx,t) = \uu_s(\xx) + \vv(\xx,t)$, the vorticity equation 
\eqref{3DNSom} can be written in equivalent form
\begin{equation}\label{3Domstrain}
  \partial_t \oomega(\xx,t) + (\vv(\xx,t),\nabla)\oomega(\xx,t) 
  - (\oomega(\xx,t),\nabla)\vv(\xx,t) \,=\, \cL \oomega(\xx,t)\,, 
\end{equation}
where $\cL$ is the linear operator defined by
\begin{equation}\label{LLdef}
  \cL\oomega \,=\, \nu\Delta \oomega - (\cM \xx,\nabla)\oomega + 
  \cM \oomega\,.
\end{equation}
As $\div\vv = 0$ and $\curl\vv = \oomega$, the Biot-Savart law
\eqref{BS3D} can be used to reconstruct the time-dependent velocity
field $\vv$ from the vorticity distribution $\oomega$. In addition to
the Laplacian, the linear operator $\cL$ includes an advection term
that depends linearly on the space variable $\xx$, and a zero order
term involving the strain matrix $\cM$ whose main effect is to amplify
the third component $\omega_3$ while attenuating $\omega_1$ and
$\omega_2$. 

The Burgers vortex is a stationary solution of \eqref{3Domstrain}
which results from the balance between the amplification of vorticity
due to stretching and the dissipation due to viscosity. In the
axisymmetric case where $\gamma_1 = \gamma_2 = -\gamma/2$ and
$\gamma_3 = \gamma > 0$, it has the explicit form
\begin{equation}\label{GGdef}
  \oomega(\xx) \,=\, \Gamma\,\frac{\gamma}{\nu} 
  \,\GG\Bigl(x\sqrt{\gamma/\nu}\Bigr)\,,
  \qquad \vv(\xx) \,=\, \Gamma\sqrt{\frac{\gamma}{\nu}}\,\vv^G\Bigl(
  x\sqrt{\gamma/\nu}\Bigr)\,,
\end{equation}
where $\Gamma$ is the total circulation and 
\begin{equation}\label{defGG}
  \GG(\xi) \,=\, \begin{pmatrix} 0 \\ 0 \\ G(\xi)\end{pmatrix}\,,
  \qquad 
  \vv^G(\xi) \,=\, \frac{1}{2\pi |\xi|^2} \Bigl(1 - e^{-|\xi|^2/4}\Bigr)
  \begin{pmatrix} -\xi_2 \\ \xi_1 \\ 0\end{pmatrix}\,.
\end{equation}
The striking similarity with the corresponding expressions
\eqref{Osdef}, \eqref{Gdef} for the Lamb-Oseen vortex is of course not
an accident. Indeed, if the pair $\oomega(\xx,t), \vv(\xx,t)$ is a
solution of \eqref{3Domstrain} that is two-dimensional in the sense
that $\partial_3 \oomega = \partial_3 \vv \equiv 0$ and $\omega_1 =
\omega_2 \equiv 0$, then the pair $\omega(x,t), u(x,t)$ defined by
\emph{Lundgren's transformation} \cite{Lu}
\begin{equation}\label{Lundgren}
  \omega(x,t) \,=\, \frac{1}{\gamma t}\,\omega_3\Bigl(\frac{x}{
  \sqrt{\gamma t}}\,,\,\frac{1}{\gamma}\log(\gamma t)\Bigr)\,,
  \qquad 
  u(x,t) \,=\, \frac{1}{\sqrt{\gamma t}}\,v\Bigl(\frac{x}{
  \sqrt{\gamma t}}\,,\,\frac{1}{\gamma}\log(\gamma t)\Bigr)\,,
\end{equation} 
satisfies the two-dimensional vorticity equation \eqref{2DNSom}.  In
other words, the two-dimensional solutions of equation
\eqref{3Domstrain}, which includes an axisymmetric linear straining
field, are in one-to-one correspondence with those of the
two-dimensional vorticity equation \eqref{2DNSom}, via a self-similar
change of variables. This observation plays a crucial role both in
Section~\ref{sec1}, where stability of the Lamb-Oseen vortex is
studied, and in Section~\ref{sec3} where the corresponding results for
the axisymmetric Burgers vortex are presented. There is however an
important difference between both situations\: although the Burgers
vortex is a two-dimensional stationary solution of \eqref{3Domstrain},
there is no reason to restrict the stability analysis to perturbations
in the same class. Quite the contrary, the Burgers vortex can be a
relevant model for tubular structures in turbulent flows only if one
can prove stability with respect to general three-dimensional
perturbations, and this is a difficult problem that has no counterpart
in the two-dimensional case, see Section~\ref{sec3} for a detailed
discussion.

In the asymmetric case where $\gamma_1 \neq \gamma_2$, Burgers
vortices still exist, but their profiles satisfy a genuinely nonlinear
equation and explicit formulas such as \eqref{GGdef}, \eqref{defGG}
are no longer available. Thus even existence of such stretched
vortices is a challenging mathematical question, which will also be
discussed in Section~\ref{sec3}. More generally, all existence, 
uniqueness, and stability results available for the axisymmetric 
Burgers vortex are expected to remain true in the asymmetric case 
too, although rigorous proofs are not always available. 

\begin{rem}
Although physical constants are useful for dimensional analysis 
and important for comparison with experiments, they often 
hinder the mathematical analysis by making formulas needlessly 
complicated. In Sections~\ref{sec1} and \ref{sec3} below, 
dimensionless variables and functions are systematically used, 
and this amounts to setting $\nu = \gamma = 1$ in all formulas. 
In particular, the total circulation of a vortex coincides with 
the circulation Reynolds number, and will be denoted by $\alpha$. 
\end{rem}

\section{Stability of Lamb-Oseen vortices}\label{sec1}

This section is devoted to the stability analysis of the family of
Lamb-Oseen vortices \eqref{Osdef}. These are self-similar solutions of
the two-dimensional vorticity equation \eqref{2DNSom}, and their
properties are most conveniently studied if the equation itself is
written in self-similar variables $\xi = x/\sqrt{t}$, $\tau =
\log(t)$ \cite{GW0}. Assuming $\nu = 1$ and setting
\begin{equation}\label{vwdef}
  \omega(x,t) \,=\, \frac{1}{t}\,w\Bigl(\frac{x}{\sqrt{t}}\,,\,
  \log(t)\Bigr)\,, \qquad 
  u(x,t) \,=\, \frac{1}{\sqrt{t}}\,v\Bigl(\frac{x}{\sqrt{t}}\,,\,
  \log(t)\Bigr)\,,
\end{equation}
one obtains for the rescaled vorticity $w(\xi,\tau)$ and the 
rescaled velocity $v(\xi,\tau)$ the following evolution 
equation
\begin{equation}\label{NS2Dw}
  \partial_\tau w(\xi,\tau) + v(\xi,\tau)\cdot\nabla w(\xi,\tau) \,=\, 
  L w(\xi,\tau)\,,
\end{equation}
where $L$ is the linear operator defined by
\begin{equation}\label{Ldef}
  L \,=\, \Delta + \frac{\xi}{2}\cdot\nabla + 1\,.
\end{equation}
The change of variables \eqref{vwdef} coincides with Lundgren's
transformation \eqref{Lundgren}, except that it is used here in the
opposite way\: starting from the two-dimensional vorticity
$\omega(x,t)$ and velocity $u(x,t)$, one obtains the rescaled
quantities $w(\xi,\tau), v(\xi,\tau)$ whose physical meaning is not
immediately obvious. In addition, the rescaled equation \eqref{NS2Dw}
looks more complicated than the original vorticity equation
\eqref{2DNSom} because the Laplace operator $\Delta$ is replaced by
the Fokker-Planck operator $L$. However, from a mathematical point of
view, the rescaled equation \eqref{NS2Dw} has several advantages which
greatly simplify the analysis. In particular, the operator $L$ has
(partially) discrete spectrum when considered in appropriate function
spaces, and that observation is crucial for the stability analysis of
the Lamb-Oseen vortex presented below. Moreover, the associated
semigroup $e^{\tau L}$ has nice confinement properties, as a
consequence of which it is possible to use compactness methods to
investigate the long-time behavior of solutions to the rescaled
vorticity equation \eqref{NS2Dw}, see \cite{GW1}.

Due to scale invariance, the Biot-Savart law \eqref{BS2D} is not 
affected by the change of variables \eqref{vwdef}. This means
that the rescaled velocity $v(\xi,\tau)$ can be reconstructed 
from the rescaled vorticity $w(\xi,\tau)$ through the formula
\begin{equation}\label{newBS2D}
  v(\cdot,\tau) \,=\, K_{2D} * w(\cdot,\tau)~, \qquad 
  \hbox{where}\quad K_{2D}(\xi) \,=\, \frac{1}{2\pi}
  \frac{\xi^\perp}{|\xi|^2}\,.
\end{equation}
By construction, for any $\alpha \in \R$, the Lamb-Oseen vortex $w =
\alpha G$, $v = \alpha v^G$ is a stationary solution of \eqref{NS2Dw}. 
The dynamical relevance of this family of equilibria is demonstrated 
by the following global convergence result.

\begin{thm}[{\bf \cite{GW1}}]\label{globconv} 
For any initial data $w_0 \in L^1(\R^2)$, the rescaled vorticity
equation \eqref{NS2Dw} has a unique global solution $w \in 
C^0([0,\infty),L^1(\R^2))$. This solution satisfies 
$\|w(\tau)\|_{L^1(\R^2)} \le \|w_0\|_{L^1(\R^2)}$ for all 
$\tau \ge 0$, and 
\begin{equation}\label{globconvest}
  \lim_{\tau \to \infty} \|w(\tau) - \alpha G\|_{L^1(\R^2)} \,=\, 
  0\,, \qquad \hbox{where}\quad \alpha \,=\, \int_{\R^2} 
  w_0(\xi)\dd \xi\,.
\end{equation}
\end{thm}

Theorem~\ref{globconv} shows that Lamb-Oseen vortices describe, to
leading order, the long-time behavior of all solutions of the
two-dimensional Navier-Stokes equations with integrable initial
vorticity and nonzero total circulation $\alpha$.  Similar conclusions
were previously obtained for small solutions \cite{GK}, and for large
solutions with small circulation \cite{Ca}. The first step in the
proof consists in showing that the original vorticity equation
\eqref{2DNSom} is globally well-posed in $L^1(\R^2)$, that the $L^1$
norm of the solutions is nonincreasing in time, and that the total
circulation $\alpha = \int_{\R^2} \omega\dd x$ is a conserved quantity
\cite{BA}.  Since the change of variables \eqref{vwdef} leaves the
$L^1$ norm invariant, the same conclusions hold for the rescaled
vorticity equation \eqref{NS2Dw} too. Then, in view of the confinement
properties of the linear semigroup $e^{\tau L}$, one can show that the
solutions of \eqref{NS2Dw} are not only bounded, but also relatively
compact in the space $L^1(\R^2)$. Finally, using appropriate Lyapunov
functions \cite{GW1} or monotonicity properties based on rearrangement
techniques \cite{GGL}, one can prove that the omega-limit set in
$L^1(\R^2)$ of any solution of \eqref{NS2Dw} is included in the family
of Lamb-Oseen vortices. As the total circulation is conserved, the
omega-limit set is in fact reduced to the singleton $\{\alpha G\}$,
which proves \eqref{globconvest}. The interested reader is referred to
\cite{GW1,GGS} for details.

The global convergence result \eqref{globconvest} is very general, 
but the proof sketched above is not constructive, and does not 
yield any estimate on the time needed to reach the asymptotic 
regime described by the Lamb-Oseen vortex. Explicit estimates 
of the convergence time can however be obtained if the vorticity 
has a definite sign \cite{GW1} or is strongly localized \cite{GR}.
In the rest of this section, emphasis is put on local stability 
results, for which explicit bounds are also available. 

\subsection{Local stability results}\label{ss1.1}

Theorem~\ref{globconv} strongly suggests, but does not really prove,
that the Lamb-Oseen vortex $\alpha G$ is a stable equilibrium of the
rescaled vorticity equation \eqref{NS2Dw} for any $\alpha \in \R$.
Stability can be established by considering solutions of the 
form $w = \alpha G + \tilde w$, $v = \alpha v^G + \tilde v$. 
The perturbations satisfy the evolution equation
\begin{equation}\label{NS2Dwper}
   \partial_\tau \tilde w + \tilde  v \cdot \nabla \tilde w \,=\, 
   (L - \alpha \Lambda)\tilde w\,, 
\end{equation}
where $L$ is given by \eqref{Ldef} and $\Lambda$ is the nonlocal 
linear operator defined by
\begin{equation}\label{Lamdef}
  \Lambda \tilde w \,=\, v^G \cdot\nabla \tilde w + 
  \tilde v\cdot\nabla G\,, \qquad \hbox{with}\quad 
  \tilde v \,=\, K_{2D}*\tilde w\,.
\end{equation} 
It is possible to prove that the perturbation equation \eqref{NS2Dwper} 
is globally well-posed in the space $L^1(\R^2)$, and that the origin
$\tilde w = 0$ is a stable equilibrium, but at this level of
generality little can be said about the long-time behavior of the
solutions. However, more precise stability results can be obtained 
if one assumes that the vorticity is sufficiently localized in space.

Given $m \in [0,\infty]$, let $\rho_m : [0,\infty) \to [1,\infty)$ be 
the weight function defined by 
\begin{equation}\label{def.rhom}
  \rho_m(r) \,=\, \begin{cases} 1 & \hbox{if} \quad m = 0\,, \\
  (1+\frac{r}{4m})^m  & \hbox{if} \quad 0 < m<\infty\,, \\
  e^{r/4} & \hbox{if} \quad m = \infty\,.\end{cases}
\end{equation}
Perturbations will be taken in the weighted $L^2$ space 
\begin{equation}\label{def.L^2(m)}
  L^2(m) \,=\, \Bigl\{w\in L^2(\R^2)~\Big|~\|w\|_{L^2(m)}^2
  = \int_{\R^2} \rho_m(|\xi|^2) |w(\xi)|^2 \dd\xi < \infty \Bigr\}\,, 
\end{equation}
which is a (real) Hilbert space equipped with the scalar product
\begin{equation}\label{L2mscalar}
  \langle w_1\,,\,w_2\rangle_{L^2(m)} \,=\, \int_{\R^2} \rho_m(|\xi|^2) 
  w_1(\xi) w_2(\xi)\dd\xi\,.
\end{equation}
Elements of $L^2(m)$ are square integrable functions with algebraic
decay at infinity if $0 < m < \infty$, and Gaussian decay if $m =
\infty$. H\"older's inequality implies that $L^2(m) \hookrightarrow 
L^1(\R^2)$ if $m > 1$. In that case, it is useful to introduce 
the closed subspace
\begin{equation}\label{def.L^2_0(m)}
  L^2_0(m) \,=\, \Bigl\{w \in L^2(m)~\Big|~\int_{\R^2} w(\xi)
  \dd\xi = 0\Bigr\}\,, 
\end{equation}
which happens to be invariant under the action of both linear 
operators $L$ and $\Lambda$. 

To study the stability of the origin $\tilde w = 0$ for the
perturbation equation \eqref{NS2Dwper}, it is useful to compute the
spectrum of the linearized operator $L - \alpha \Lambda$ in the
(complexified) Hilbert space $L^2(m)$. In the simple case where
$\alpha = 0$, the spectrum is explicitly known\:

\begin{prop}[\cite{GW0}]\label{Lspectrum}
For any $m \in [0,\infty]$, the spectrum of the linear operator
\eqref{Ldef} in the weighted space $L^2(m)$ defined by 
\eqref{def.rhom}, \eqref{def.L^2(m)} is 
\begin{equation}\label{Lspec}
  \sigma_m(L) \,=\, \Bigl\{\lambda \in \C\,\Big|\, \Re(\lambda) 
  \le \frac{1}{2} - \frac{m}{2}\Bigr\} \,\cup\, \Bigl\{-\frac{k}{2} 
  \,\Big|\, k \in \N\Bigr\}\,.
\end{equation}
Moreover, if $m > k+1$ for some $k \in \N$, then $\lambda_k = -k/2$ is 
an isolated eigenvalue of $L$, with (algebraic and geometric) 
multiplicity $k+1$.
\end{prop}

It follows in particular from Proposition~\ref{Lspectrum} that $L$ 
has purely discrete spectrum in $L^2(m)$ when $m = \infty$. This is
easily understood if one observes that $\rho_\infty(|\xi|^2) = e^{|\xi|^2/4} 
= (4\pi)^{-1}G(\xi)^{-1}$, and that 
\begin{equation}\label{Lconj}
  G^{-1/2} \,L ~G^{1/2} \,=\, \Delta \,-\, \frac{|\xi|^2}{16} \,+\, 
  \frac12\,.
\end{equation}
The formal relation \eqref{Lconj} implies that the operator $L$ in
$L^2(\infty)$ is unitarily equivalent to the harmonic oscillator
$\Delta - |\xi|^2/16 +1/2$ in $L^2(\R^2)$, the spectrum of which is
the sequence $(\lambda_k)_{k \in \N}$, where $\lambda_k = -k/2$ has
multiplicity $k+1$. If $m < \infty$, the discrete part of the spectrum
persists because the corresponding eigenfunctions decay rapidly at
infinity. In addition, any $\lambda \in \C$ such that $\Re(\lambda) <
(1-m)/2$ is an eigenvalue of $L$ in $L^2(m)$ with infinite
multiplicity \cite{GW0}, hence the spectrum $\sigma_m(L)$ also
includes the closed half-plane $H_m = \{\lambda \in \C\,|\,
\Re(\lambda) \le (1-m)/2\}$.

In the more interesting case where $\alpha \neq 0$, the spectrum of
$L-\alpha\Lambda$ in $L^2(m)$ cannot be computed explicitly. However, 
upper bounds on the real part of the spectrum are sufficient for 
the stability analysis, and such estimates can be obtained by 
combining the following three observations. 

\medskip\noindent{\bf Observation 1\:} The operator $\Lambda$ is a
{\em relatively compact} perturbation of $L$ in $L^2(m)$, for any $m
\in [0,\infty]$.  This is intuitively obvious, because $\Lambda$ is a
first-order differential operator whose coefficients decay to zero at
infinity, whereas $L$ involves in particular the Laplace operator
$\Delta$. By Weyl's theorem, the {\em essential spectrum} \cite{He} 
of $L-\alpha\Lambda$ in $L^2(m)$ does not depend on $\alpha$, hence 
coincides with the closed half-plane $H_m = \{\lambda \in \C\,|\, 
\Re(\lambda) \le (1-m)/2\}$ by Proposition~\ref{Lspectrum}. It 
thus remains to locate isolated eigenvalues of $L-\alpha\Lambda$ 
outside $H_m$. 

\medskip\noindent{\bf Observation 2\:} The isolated eigenvalues of
$L-\alpha\Lambda$ in $L^2(m)$ do not depend on $m$. Indeed, if $w \in
L^2(m)$ satisfies $(L-\alpha\Lambda)w = \lambda w$ for some $\lambda
\in \C \setminus H_m$, one can show that $w$ decays sufficiently fast
at infinity so that $w \in L^2(\infty)$ \cite{GW1}. This means that
isolated eigenvalues of $L-\alpha\Lambda$ can be located by
considering the particular case $m = \infty$, where the spectrum 
is fully discrete and consists of a sequence of eigenvalues 
$(\lambda_k(\alpha))_{k \in \N}$ with $\Re(\lambda_k(\alpha)) \to 
-\infty$ as $k \to \infty$. 

\medskip\noindent{\bf Observation 3\:} The operator $\Lambda$ 
is {\em skew-symmetric} in $L^2(\infty)$, namely  
\begin{equation}\label{Lamasym}
  \langle \Lambda w_1\,,\,w_2\rangle + \langle w_1\,,\,\Lambda 
  w_2\rangle \,=\, 0\,, \qquad \hbox{for all} \quad w_1,w_2 \in 
  D(\Lambda) \subset L^2(\infty)\,,
\end{equation}
where $D(\Lambda) \subset L^2(\infty)$ is the (maximal) domain of the
operator $\Lambda$, and $\langle\cdot\,,\cdot\rangle$ denotes the
scalar product in $L^2(\infty)$, which (up to an irrelevant factor)
can be written in the form
\begin{equation}\label{Xprod}
  \langle w_1 \,,\,w_2\rangle \,=\, \int_{\R^2} G(\xi)^{-1}
  w_1(\xi)w_2(\xi)\dd\xi\,.
\end{equation}
To prove \eqref{Lamasym} one decomposes $\Lambda = \Lambda_1 +
\Lambda_2$, where $\Lambda_1 w = v^G \cdot \nabla w$ and $\Lambda_2 w
= (K_{2D}*w)\cdot \nabla G$. If $w_1, w_2 \in L^2(\infty)$ belong to
the domain of $\Lambda$, then
\begin{align*}
  \langle \Lambda_1 w_1,w_2\rangle +  \langle w_1, \Lambda_1 w_2
  \rangle \,&=\, \int_{\R^2} G^{-1}\Bigl(w_2 \,v^G \cdot \nabla w_1
  + w_1\,v^G\cdot \nabla w_2\Bigr)\dd\xi \cr
  \,&=\, \int_{\R^2} G^{-1} \,v^G\cdot \nabla(w_1 w_2)\dd\xi 
  \,=\, 0\,, 
\end{align*}
because the vector field $G^{-1}v^G$ is divergence-free. Moreover
using the identity $\nabla G = -\frac12\xi G$ and the Biot-Savart 
law \eqref{newBS2D}, one obtains
\begin{align*}
  &\langle \Lambda_2 w_1,w_2\rangle +  \langle w_1, \Lambda_2 w_2
  \rangle \,=\, -\frac12\int_{\R^2} \Bigl((\xi\cdot v_1) w_2 
  + (\xi\cdot v_2) w_1\Bigr)\dd\xi \cr
  \,&=\, -\frac{1}{4\pi}\int_{\R^2}\int_{\R^2}
  \left\{\xi\cdot \frac{(\xi-\eta)^\perp}{|\xi-\eta|^2} +
  \eta \cdot \frac{(\eta-\xi)^\perp}{|\xi-\eta|^2}\right\}
  w_1(\eta) w_2(\xi)\dd\eta\dd\xi \,=\, 0\,,
\end{align*}
because the last integrand vanishes identically. This proves
\eqref{Lamasym}. One can also show that the operator $\Lambda$ 
is not only skew-symmetric, but also skew-adjoint in $L^2(\infty)$, 
see \cite{M3}.

The observations above lead to the following {\em spectral stability
result} for the Lamb-Oseen vortex in the space $L^2(m)$. 

\begin{prop}[\cite{GW1}]\label{specstab}
For any $\alpha \in \R$ and any $m \in [1,\infty]$, the spectrum 
of the linearized operator $L - \alpha\Lambda$ in the space $L^2(m)$ 
satisfies
\begin{equation}\label{sigma0}
  \sigma_m(L-\alpha\Lambda) \,\subset\, \Bigl\{\lambda \in \C
  \,\Big|\, \Re(\lambda) \le 0\Bigr\}\,.
\end{equation}
Moreover, if $m \ge 2$, then 
\begin{equation}\label{sigma1}
  \sigma_m(L-\alpha\Lambda) \,\subset\, \{0\} \,\cup\, \Bigl\{\lambda 
  \in \C\,\Big|\, \Re(\lambda) \le -\frac12\Bigr\}\,.
\end{equation}
Finally, if $m \ge 3$, then 
\begin{equation}\label{sigma2}  
  \sigma_m(L-\alpha\Lambda) \,\subset\, \{0\} \,\cup\, 
  \Bigl\{-\frac12\Bigr\} \,\cup\, \Bigl\{\lambda \in \C \,\Big|\, 
  \Re(\lambda) \le -1\Bigr\}\,.
\end{equation}
\end{prop}

\noindent{\bf Proof.}
As before let $H_m = \{\lambda \in \C\,|\, \Re(\lambda) \le (1-m)/2\}$. 
By Observation 1 above, if $m \ge 1$, the essential spectrum of 
$L-\alpha\Lambda$ is included in the half-space $H_1$. Assume 
that $\lambda \in \C \setminus H_1$ is an isolated eigenvalue 
of $L-\alpha\Lambda$, and let $w \in L^2(m)$ be a nontrivial  
eigenfunction associated with $\lambda$. Then $w \in L^2(\infty)$ 
by Observation 2, and using Observation 3 one finds
\begin{equation}\label{dissip}
  \Re(\lambda) \langle w,w\rangle \,=\, \Re \langle 
  (L-\alpha\Lambda)w,w\rangle \,=\, \langle L w,w\rangle 
  \,\le\, 0\,,
\end{equation}
because $L$ is a nonpositive self-adjoint operator in $L^2(\infty)$
and $\Lambda$ is skew-symmetric. This contradicts the assumption that
$\Re(\lambda) > 0$, hence the whole spectrum of $L-\alpha\Lambda$ in
$L^2(m)$ is contained in the half-space $H_1$, as asserted in
\eqref{sigma0}.

As $LG = \Lambda G = 0$, it is clear that $0$ is an eigenvalue of
$L-\alpha\Lambda$ for any $m \ge 0$ and any $\alpha \in \R$.  If $m >
1$, one can write $L^2(m) = \R G \oplus L^2_0(m)$, where $L^2_0(m)$
is the hyperplane defined in \eqref{def.L^2_0(m)}, and this 
decomposition is left invariant by both operators $L$ and 
$\Lambda$. Now, the same argument as above shows that, if $m \ge
2$, the spectrum of the operator $L-\alpha\Lambda$ acting on 
$L^2_0(m)$ is contained in the half-plane $H_2$, because $L \le -1/2$ 
on $L^2_0(\infty)$. This proves \eqref{sigma1}.

Finally, it is easy to verify that $L(\partial_i G) = -\frac12 
\partial_i G$ for $i = 1,2$, and differentiating the identity 
$v^G \cdot\nabla G = 0$ one finds that $\Lambda (\partial_i G) = 0$ 
for $i = 1,2$. This means that $-1/2$ is an eigenvalue of 
$L-\alpha\Lambda$ for any $m \ge 0$ and any $\alpha \in \R$.
As above, if $m > 2$, one has the invariant decomposition
\[
  L^2(m) \,=\,  \{\alpha G\,|\, \alpha \in \R\} \oplus 
  \{\beta_1\partial_1 G + \beta_2\partial_2 G\,|\, \beta_1,\beta_2 
  \in \R\} \oplus L^2_{00}(m)~,
\]
where
\begin{equation}\label{def.L^2_00(m)}
  L^2_{00}(m) \,=\, \Bigl\{w \in L^2_0(m)~\Big|~\int_{\R^2} \xi_i w(\xi)
  \dd\xi = 0 \hbox{ for }i = 1,2\Bigr\}\,. 
\end{equation}
As $L \le -1$ on $L^2_{00}(\infty)$, the same argument shows that 
the spectrum of the operator $L-\alpha\Lambda$ acting on $L^2_{00}(m)$ 
is contained in the half-plane $H_3$, if $m \ge 3$. This proves 
\eqref{sigma2}. 
\QED

\medskip 
The linear operator $L-\alpha\Lambda$ is the generator of a
strongly continuous semigroup in the space $L^2(m)$ for any $\alpha
\in \R$ and any $m \in [0,\infty]$ \cite{GW0}. The following {\em linear
  stability result} is a natural consequence of Proposition~\ref{specstab} 
and its proof.

\begin{prop}[\cite{GW1}]\label{linstab}
For any $\alpha \in \R$ and any $m > 1$, there exists a 
positive constant $C$ such that
\begin{equation}\label{linstab0}
  \|e^{\tau(L-\alpha\Lambda)}\|_{L^2(m) \to L^2(m)} \,\le\, C\,,\qquad 
  \hbox{for all } \tau \ge 0\,.
\end{equation}
Moreover, if $m > 2$, then 
\begin{equation}\label{linstab1}
  \|e^{\tau(L-\alpha\Lambda)}\|_{L^2_0(m) \to L^2_0(m)} \,\le\, C\,e^{-\tau/2}\,,
  \qquad \hbox{for all } \tau \ge 0\,.
\end{equation}
Finally, if $m > 3$, then
\begin{equation}\label{linstab2}
  \|e^{\tau(L-\alpha\Lambda)}\|_{L^2_{00}(m) \to L^2_{00}(m)} \,\le\, 
  C\,e^{-\tau}\,,\qquad \hbox{for all } \tau \ge 0\,.
\end{equation}
\end{prop}

When studying the stability of the Lamb-Oseen vortex, there is no loss
of generality in considering perturbations with zero total
circulation. Indeed, if $w = \alpha G + \tilde w$ for some $\tilde w
\in L^2(m)$ with $m > 1$, then defining $\tilde \alpha = \int_{\R^2}
\tilde w(\xi)\dd\xi$ one can write $w = (\alpha+\tilde \alpha)G +
(\tilde w - \tilde \alpha G)$, where by construction $\tilde w -
\tilde \alpha G \in L^2_0(m)$. Thus perturbations with nonzero
circulation of the vortex $\alpha G$ can be considered as
perturbations with zero circulation of the modified vortex $(\alpha +
\tilde \alpha)G$.  As the total circulation is a conserved quantity,
the subspace $L^2_0(m)$ is invariant under the evolution defined by
the full perturbation equation \eqref{NS2Dwper}. By
Proposition~\ref{linstab}, the linear semigroup
$e^{\tau(L-\alpha\Lambda)}$ is exponentially decaying in $L^2_0(m)$ if
$m > 2$, and using that information it is routine to deduce the
following {\em asymptotic stability result}, which is the main outcome
of this section.

\begin{prop}[\cite{GW1}]\label{nonlinstab}
Fix $\alpha \in \R$ and $m \in (2,\infty]$. There exist positive 
constants $\epsilon$ and $C$ such that, for all $\tilde w_0 \in 
L^2_0(m)$ satisfying $\|\tilde w_0\|_{L^2(m)} \le \epsilon$, the 
rescaled vorticity equation \eqref{NS2Dw} has a unique global 
solution $w \in C^0([0,\infty),L^2(m))$ with initial data
$w_0 = \alpha G + \tilde w_0$. Moreover, the following estimate
holds
\begin{equation}\label{nonlindecay}
  \|w(\tau) - \alpha G\|_{L^2(m)} \,\le\, C \|w_0 - \alpha G\|_{L^2(m)}
  \,e^{-\tau/2}~, \qquad \tau \ge 0~.
\end{equation}
\end{prop}

If $m > 2$, the codimension $3$ subspace $L^2_{00}(m)$ is also 
invariant under the evolution defined by the full perturbation 
equation \eqref{NS2Dw}. As a consequence, if $\tilde w_0 \in 
L^2_{00}(m)$, the solution of \eqref{NS2Dw} given by 
Proposition~\ref{nonlinstab} satisfies $w(\tau) - \alpha G 
\in L^2_{00}(m)$ for all $\tau \ge 0$. If $m > 3$, one can then 
use \eqref{linstab2} to conclude that
\[
  \|w(\tau) - \alpha G\|_{L^2(m)} \,\le\, C \|w_0 - \alpha G\|_{L^2(m)}
  \,e^{-\tau}~, \qquad \tau \ge 0~.
\]
As is shown in \cite{GW1,Ga1}, if $\alpha \neq 0$, the assumption
that $w_0$ has vanishing first order moments does not really 
restrict the generality, because this condition can always be met 
by a suitable translation of the initial data. 

\begin{rem}\label{uniformstab}
If $m = \infty$, one can show that the Lamb-Oseen vortex $\alpha G$ is
uniformly stable for all $\alpha \in \R$ in the sense that the
constants $\epsilon$ and $C$ in Proposition~\ref{nonlinstab} do not
depend on $\alpha$ \cite{Ga1}. This is in sharp contrast with what
happens for shear flows, such as the Poiseuille flow in a cylindrical
pipe or the Couette-Taylor flow between two rotating cylinders. In
such examples, the laminar stationary flow undergoes an instability,
of spectral or pseudospectral nature, when the Reynolds number is
sufficiently large. In contrast, a fast rotation has rather a
stabilizing effect on vortices, as the analysis below reveals.
\end{rem}

\subsection{Large Reynolds number asymptotics}\label{ss1.2}

Proposition~\ref{specstab} above gives uniform estimates on the
spectrum of the linearized operator $L - \alpha\Lambda$, which are
sufficient to prove stability of the Lamb-Oseen vortex for all values
of the circulation parameter $\alpha \in \R$. However, such estimates
do not describe how the spectrum changes as the circulation parameter
varies. The most relevant regime for turbulent flows is of course the
high Reynolds number limit where $|\alpha| \to \infty$, which deserves
a special consideration. As the essential spectrum of
$L-\alpha\Lambda$ in the space $L^2(m)$ does not depend on $m$, it is
most convenient to work in the limiting space $X = L^2(\infty)$,
equipped with the scalar product \eqref{Xprod}. In that space, as was
already mentioned, the spectrum of $L-\alpha\Lambda$ is discrete, and
consists of a sequence of eigenvalues $(\lambda_k(\alpha))_{k \in \N}$
with $\Re(\lambda_k(\alpha)) \to -\infty$ as $k \to \infty$. It
follows from Proposition~\ref{Lspectrum} that $\lambda_k(0) = -k/2$,
for any $k \in \N$, and the goal of this section is to investigate the
behavior of the real part of $\lambda_k(\alpha)$ as $|\alpha| \to
\infty$.

The starting point of the analysis is the determination of the kernel
of the skew-symmetric operator $\Lambda$. Let $X_0 \subset X$ denote
the closed subspace containing all radially symmetric functions. If
$w \in X_0$, the associated velocity field $v = K_{2D}*w$ satisfies
$\xi\cdot v(\xi) = 0$, and it follows that $\Lambda w = 0$, hence $X_0
\subset \ker(\Lambda)$. On the other hand, it was already observed
that $\Lambda (\partial_i G) = 0$ for $i = 1,2$. The following result
asserts that the kernel of $\Lambda$ does not contain any more
elements\:

\begin{lem}[\cite{M3}]\label{kerLam} $\ker(\Lambda) = X_0 
\oplus \{\beta_1\partial_1 G + \beta_2 \partial_2 G\,|\, \beta_1,
\beta_2 \in \R\}$.
\end{lem}

In view of Lemma~\ref{kerLam}, the subspace $\ker(\Lambda) \subset X$
is invariant under the action of both operators $L$ and $\Lambda$, and
the orthogonal complement $\ker(\Lambda)^\perp$ is invariant too
because $L$ is self-adjoint and $\Lambda$ is skew-adjoint. Inside
$\ker(\Lambda)$, the spectrum of $L-\alpha\Lambda \equiv L$ does not
depend on the circulation parameter $\alpha$, and consists of all
negative integers in addition to the double eigenvalue $-1/2$. In
fact, for any $n \in \N$, the eigenfunction corresponding to the
eigenvalue $-n$ is the radially symmetric Hermite function $\Delta^n
G$. The only difficult task is therefore to study the spectrum of
$L_\perp- \alpha\Lambda_\perp$, which is defined as the restriction of
$L-\alpha\Lambda$ to the orthogonal complement
$\ker(\Lambda)^\perp$. That spectrum does depend in a nontrivial way
upon the parameter $\alpha$. It happens that the real parts of all
eigenvalues converge to $-\infty$ as $|\alpha| \to \infty$, which is
of course compatible with the uniform bounds given by
Proposition~\ref{specstab}. This phenomenon illustrates the {\em
  stabilizing effect} of fast rotation on Lamb-Oseen vortices.

Two natural quantities can be introduced to accurately measure 
the effect of fast rotation. For any $\alpha \in \R$, one can define 
the {\em spectral lower bound}
\[ 
  \Sigma(\alpha) \,=\, \inf\Bigl\{\Re(z)\,\Big|\, z \in 
  \mathrm{spec}(-\cL_\perp + \alpha\Lambda_\perp)\Bigr\}~,
\]
or the {\em pseudospectral bound} 
\[
  \Psi(\alpha) \,=\, \Bigl(\sup_{\lambda \in \R}\|(\cL_\perp
  -\alpha\Lambda_\perp-i\lambda)^{-1}\|_{X \to X}\Bigr)^{-1}~.
\]
In the definition of $\Sigma(\alpha)$, the sign of the linearized
operator has been changed to obtain a positive quantity. Although
the spectral and pseudospectral bounds are of rather different 
nature, there is a simple one-sided relation between them\:

\begin{lem}\label{SigPsi}
For any $\alpha \in \R$ one has $\Sigma(\alpha) \ge \Psi(\alpha)
\ge 1$. 
\end{lem}

\noindent{\bf Proof.}
Fix $\alpha \in \R$. By Lemma~\ref{kerLam}, one has $\ker(\Lambda)^\perp 
\subset L^2_{00}(m)$, hence $\Sigma(\alpha) \ge 1$ by \eqref{sigma2}. 
On the other hand, if $(L - \alpha\Lambda + \lambda)w = 0$ for some 
$\lambda \in \C$ and some $w \in \ker(\Lambda)^\perp$ such 
that $\langle w,w\rangle = 1$, then $(L - \alpha\Lambda + 
i\Im(\lambda))w = -\Re(\lambda)w$, hence
\[
  \Re(\lambda) \,\ge\, \|(L_\perp - \alpha\Lambda_\perp +
  i\Im(\lambda))^{-1}\|^{-1} \,\ge\, \Psi(\alpha)~.
\]
This proves that $\Sigma(\alpha) \ge \Psi(\alpha)$. Finally, the proof
of Proposition~\ref{specstab} shows that the operator $L_\perp -
\alpha\Lambda_\perp +1$ is $m$-dissipative \cite{Ka}. This in
particular implies that $\|(L_\perp - \alpha\Lambda_\perp-i\lambda)^{-1}\| 
\le 1$ for all $\lambda \in \R$, hence $\Psi(\alpha) \ge 1$.
\QED

\medskip
The stabilizing effect in the large Reynolds number limit is 
qualitatively illustrated by the following result\:

\begin{prop}[\cite{M3}]\label{stabeff}
One has $\Psi(\alpha) \to \infty$ and $\Sigma(\alpha) \to \infty$ 
as $|\alpha| \to \infty$. 
\end{prop}

The proof given in \cite{M3} actually shows that $\Sigma(\alpha) \to
\infty$ as $|\alpha| \to \infty$, but can be easily modified to yield
the stronger conclusion that $\Psi(\alpha) \to \infty$.  For the
stability analysis of the Lamb-Oseen vortex $\alpha G$, the divergence
of the spectral bound means that the decay rate in time of
perturbations in $\ker(\Lambda)^\perp$ becomes arbitrarily large as
$|\alpha| \to \infty$. On the other hand, using the divergence of the
pseudospectral bound, one can show that the basin of attraction of the
Lamb-Oseen vortex, in the weighted space $L^2(\infty)$, becomes
arbitrarily large as $|\alpha| \to \infty$. It should be emphasized, 
however, that the argument used in \cite{M3} is nonconstructive and 
does not provide any explicit estimate on the quantities $\Psi(\alpha)$ 
or $\Sigma(\alpha)$ for large $|\alpha|$.  

In fact, there are good reasons to conjecture that $\Sigma(\alpha) =
\cO(|\alpha|^{1/2})$ and $\Psi(\alpha) = \cO(|\alpha|^{1/3})$ as
$|\alpha| \to \infty$. First of all, extensive numerical calculations
performed by Prochazka and Pullin \cite{PP1,PP2} indicate that
$\Sigma(\alpha) = \cO(|\alpha|^{1/2})$ as $|\alpha| \to \infty$. Next,
the conjecture is clearly supported by rigorous analytical results on
model problems \cite{GGN}. In particular, for the simplified linear
operator $L - \alpha\Lambda_1$ where the nonlocal part $\Lambda_2$ has
been omitted, it can be proved that $\Psi(\alpha) = \cO(|\alpha|^{1/3})$ 
as $|\alpha| \to \infty$ \cite{De1}. The same result holds for the
full linearized operator $L-\alpha\Lambda$ restricted to a 
smaller subspace than $\ker(\Lambda)^\perp$, where a finite number
of Fourier modes with respect to the angular variable in polar
coordinates have been removed \cite{De2}. The general case is 
still under investigation \cite{GG}. 

Assuming that the conjecture above is true, it is worth noting that
the pseudospectral bound $\Psi(\alpha)$ and the spectral bound
$\Sigma(\alpha)$ have different growth rates as $|\alpha| \to \infty$.
This reflects the fact that the linearized operator $L-\alpha\Lambda$
becomes highly non-selfadjoint in the fast rotation limit. Indeed, for
selfadjoint or normal operators, it is easy to verify that the
spectral and pseudospectral bounds always coincide.

\subsection{Lamb-Oseen vortices in exterior domains }\label{ss1.3}

As was already mentioned, the Lamb-Oseen vortex plays a double role in
the dynamics of the Navier-Stokes equations in the whole space
$\R^2$\: it is the unique solution of the system when the initial
vorticity is a Dirac measure, and it describes the long-time
asymptotics of all solutions for which the vorticity distribution is
integrable and has nonzero total circulation. The proofs given in
\cite{GW1} demonstrate that both properties are closely related, 
due to scale invariance. Now, if the fluid is contained in a
two-dimensional domain $\Omega \subset \R^2$, and satisfies (for
instance) no-slip boundary conditions on $\partial\Omega$, scale
invariance is broken and there is no simple relation anymore between
the Cauchy problem for singular initial data and the long-time
asymptotics of general solutions. Both questions are interesting and,
at present time, largely open. In this section, the relatively simple
case of a two-dimensional exterior domain is considered, where a few
results concerning the long-time behavior of solutions with nonzero
circulation at infinity have been obtained recently.

Let $\Omega \subset \R^2$ be a smooth exterior domain, namely an
unbounded connected open set with a smooth compact boundary $\partial
\Omega$. The Navier-Stokes equations in $\Omega$ with no-slip boundary
conditions can be written in the following form\:
\begin{equation}\label{NSu}
  \left\{\begin{array}{llll}
  \partial_t u + (u\cdot\nabla)u \,=\, \Delta u - \nabla p~, 
  \quad \div u \,=\, 0~, & \quad & \hbox{for~} x \in \Omega\,, \quad
  &t > 0\,, \\
  u(x,t) \,=\, 0~, & & \hbox{for~} x \in \partial\Omega\,,
  &t > 0\,, \\
  u(x,0) \,=\, u_0(x)~, & & \hbox{for~} x \in \Omega\,,&
  \end{array}\right.
\end{equation}
where $p$ denotes the ratio of the pressure to the fluid density. The
vorticity $\omega = \partial_1 u_2 - \partial_2 u_1$ still satisfies
the simple evolution equation \eqref{2DNSom} (with $\nu =1$), but the
assumption that $u = 0$ on $\partial\Omega$ translates into a
nonlinear, nonlocal boundary condition for $\omega$, which is very
difficult to handle. So, whenever possible, it is preferable to work
directly with the velocity formulation \eqref{NSu}.

If the initial velocity $u_0$ belongs to the energy space
\[
  L^2_\sigma(\Omega) \,=\, \Bigl\{u \in L^2(\Omega)^2 \,\Big|\, 
  \div u = 0 \hbox{ in }\Omega\,,~ u\cdot n = 0 \hbox{ on }
  \partial\Omega\Bigr\}\,, 
\]
where $n$ denotes the unit normal on $\partial\Omega$, it is well
known that system \eqref{NSu} has a unique global solution satisfying
the energy identity
\[
  \frac12 \|u(\cdot,t)\|_{L^2(\Omega)}^2 + \int_0^t \|\nabla u(\cdot,
  s)\|_{L^2(\Omega)}^2\dd s \,=\, \frac12 \|u_0\|_{L^2(\Omega)}^2\,,
  \qquad t \ge 0\,.
\]
That solution converges to zero in $L^2_\sigma(\Omega)$ as 
$t \to \infty$ \cite{Mas}, which means that the long-time behavior 
of all finite energy solutions is trivial. However, as in the
whole plane $\R^2$, one can consider flows with nonzero circulation 
at infinity\:
\[
  \alpha \,=\, \lim_{R \to \infty} \oint_{|x| = R} \left(u_1 \dd x_1 
  + u_2 \dd x_2\right) \,\neq\, 0\,,
\]
in which case the kinetic energy is necessarily infinite and
the long-time behavior may be nontrivial. 

To construct such solutions, it is convenient to introduce a smooth
cut-off function $\chi : \R^2 \to [0,1]$ such that $\chi$ vanishes in
a neighborhood of $\R^2 \setminus \Omega$ and $\chi(x) = 1$ whenever
$|x|$ is sufficiently large. For technical reasons, one also assumes
that $\chi$ is radially symmetric and nondecreasing along rays. The
{\em truncated Oseen vortex}
\begin{equation}\label{uchidef}
  u^\chi(x,t) \,=\, \frac{1}{2\pi} \,\frac{x^\perp}{|x|^2} 
  \Bigl(1 - e^{-\frac{|x|^2}{4(1+t)}}\Bigr)\chi(x)\,, \qquad x \in \R^2, 
  \quad t \ge 0\,,
\end{equation}
is a divergence-free velocity field which vanishes identically 
in a neighborhood of $\R^2 \setminus \Omega$ and coincides with
the Lamb-Oseen vortex (with unit circulation) far away from 
the origin. In particular $u^\chi \notin L^2(\Omega)$. 
The corresponding vorticity distribution $\omega^\chi = 
\partial_1 u^\chi_2 - \partial_2 u^\chi_1$ reads\:
\begin{equation}\label{omegachi}
  \omega^\chi(x,t) \,=\, \frac{1}{4\pi(1+t)}\,e^{-\frac{|x|^2}{4(1+t)}}
  \,\chi(x) + \frac{1}{2\pi}\,\frac{1}{|x|^2}\Bigl(1 - 
  e^{-\frac{|x|^2}{4(1+t)}}\Bigr)x\cdot\nabla\chi(x)\,,
\end{equation}
and satisfies $\int_\Omega \omega^\chi(x,t)\dd x = 1$ for all $t \ge
0$. Of course the velocity field $u^\chi$ is not an exact solution of
the Navier-Stokes equations \eqref{NSu} (unless $\Omega = \R^2$ and
$\chi \equiv 1$), but the following result shows that it is 
a globally stable asymptotic solution. 

\begin{thm}[\cite{GM2,M4}]\label{OseenExt}
Fix $q \in (1,2]$, and let $\mu = 1/q - 1/2$. There exists a 
constant $\epsilon > 0$ such that, for all initial data of the 
form $u_0 = \alpha u^\chi(\cdot,0) + v_0$ with $|\alpha| \le \epsilon$ 
and $v_0 \in L^2_\sigma(\Omega) \cap L^q(\Omega)^2$, the Navier-Stokes 
equations \eqref{NSu} have a unique global solution which 
satisfies
\begin{equation}\label{NSuconv}
  \|u(\cdot,t) -\alpha u^\chi(\cdot,t)\|_{L^2(\Omega)} + 
  t^{1/2} \|\nabla u(\cdot,t) -\alpha\nabla u^\chi(\cdot,t)
  \|_{L^2(\Omega)} \,=\, \cO(t^{-\mu})~,
\end{equation}
as $t \to +\infty$. Moreover, if $q = 2$, then $\|u(\cdot,t) -\alpha 
u^\chi(\cdot,t)\|_{L^2(\Omega)} \to 0$ as $t \to +\infty$. 
\end{thm}

Several comments are in order. Existence and uniqueness of global
solutions to the Navier-Stokes equations \eqref{NSu} for a class of
infinite-energy initial data including those considered in
Theorem~\ref{OseenExt} were established by Kozono and Yamazaki in
\cite{KY}. The novelty here is the description of the long-time
asymptotics for a specific family of solutions, corresponding to
spatially localized perturbations of the truncated Oseen
vortex. Theorem~\ref{OseenExt} is a global stability result, in the
sense that arbitrary large perturbations $v_0$ of the vortex can be
considered. There is, however, a limitation on the size of the
circulation parameter $\alpha$, which is probably of technical
nature. To remove that restriction it seems rather natural to use the
vorticity formulation and the nice properties of the linearized
operator established in Section~\ref{ss1.1}, but this is difficult in
the present case because the boundary condition for the vorticity is
very awkward. The parameter $q$ in Theorem~\ref{OseenExt} measures the
spatial decay of the initial perturbations $v_0$ to the Oseen vortex,
and is directly related to the decay rate in time (called $\mu$) of
the corresponding solutions. If $q < 2$, then $\mu > 0$ and it is
shown in \cite{GM2} that the constant $\epsilon$ depends only on $q$,
and not on the domain $\Omega$. The limiting case where $q = 2$ was
treated in \cite{IKL,M4}.

The main original ingredient in the proof of Theorem~\ref{OseenExt} 
is a logarithmic energy estimate that is worth discussing briefly.
For solutions of the Navier-Stokes equations \eqref{NSu} of 
the form $u(x,t) = \alpha u^\chi(x,t) + v(x,t)$, the perturbation
$v$ satisfies
\begin{equation}\label{NSv}
  \partial_t v + \alpha (u^\chi,\nabla)v + \alpha (v,\nabla)u^\chi
  + (v,\nabla)v \,=\, \Delta v + \alpha R^\chi - \nabla q\,, 
  \qquad \div v \,=\, 0\,,
\end{equation}
where the source term $R^\chi = \Delta u^\chi - \partial_t u^\chi$ 
measures by how much the truncated vortex $u^\chi$ fails to 
be an exact solution of \eqref{NSu}. Taking into account the 
uniform bounds
\[
  \|\nabla u^\chi(\cdot,t)\|_{L^\infty(\R^2)} \,\le\, \frac{b}{1+t}\,,\qquad
  \|R^\chi(\cdot,t)\|_{L^2(\R^2)} \,\le\, \frac{\kappa}{1+t}\,,
\]
which hold for some positive constants $b,\kappa$ depending 
only on the cut-off $\chi$, a standard energy estimate yields
the differential inequality
\begin{equation}\label{roughenergy}
  \frac12 \frac{\D}{\D t} \|v(t)\|_{L^2(\Omega)}^2 + \|\nabla v(t)
  \|_{L^2(\Omega)}^2 \,\le\, \frac{b|\alpha|}{1+t}\,\|v(t)\|_{L^2(\Omega)}^2
  + \frac{\kappa|\alpha|}{1+t}\,\,\|v(t)\|_{L^2(\Omega)}\,,
\end{equation}
which predicts a polynomial growth of the $L^2$ norm 
$\|v(t)\|_{L^2(\Omega)}$. This naive estimate can be substantially
improved if one observes that the truncated Oseen vortex $u^\chi$ decays
like $|x|^{-1}$ as $|x| \to \infty$, and thus nearly belongs to 
the energy space. The optimal result is\:

\begin{prop}\label{logbound} There exists a constant $K > 0$ 
such that, for any $\alpha \in \R$ and any $v_0 \in L^2_\sigma(\Omega)$, 
the solution of \eqref{NSv} with initial data $v_0$ satisfies, 
for all $t \ge 1$,
\begin{equation}\label{vlog}
  \|v(t)\|_{L^2(\Omega)}^2 + \int_0^t \|\nabla v(s)\|_{L^2(\Omega)}^2\dd s 
  \,\le\, K \Bigl(\|v_0\|_{L^2(\Omega)}^2 + \alpha^2 \log(1+t) +
  \alpha^2\log(2+|\alpha|)\Bigr)\,.
\end{equation}
\end{prop}

In the proof of Theorem~\ref{OseenExt}, the logarithmic bound
\eqref{vlog} is combined with standard energy estimates for a
fractional primitive of the velocity field to prove that $v(\cdot,t)$
converges to zero in $L^2(\Omega)$ as $t \to \infty$, see
\cite{GM2,Ga2}. The optimal decay rate in \eqref{NSuconv} is then
obtained by a direct study of small solutions to the perturbation
equation \eqref{NSv}.

Although Theorem~\ref{OseenExt} is established using the velocity
formulation of the Navier-Stokes system, it is instructive to see what
it implies for the vorticity distribution $\omega$.  Assume for
instance that the initial vorticity $\omega_0 = \partial_1 (u_0)_2 
-\partial_2 (u_0)_1$ is sufficiently localized so that
\[
  \int_{\Omega} (1 + |x|^2)^m |\omega_0 (x)|^2 \dd x \,<\,\infty\,, 
\]
for some $m > 1$. By H\"older's inequality this implies that $\omega_0
\in L^1(\Omega)$. Then denoting $v_0 = u_0 - \alpha u^\chi(\cdot,0)$
where $\alpha = \int_{\Omega} \omega_0 (x)\dd x$, it follows that $v_0
\in L^2_\sigma(\Omega) \cap L^q(\Omega)^2$ for any $q \in (1,2)$ such
that $q > 2/m$ \cite{GW0}. In particular, if $|\alpha|\le \epsilon$,
the conclusion of Theorem \ref{OseenExt} holds. Moreover, the
vorticity satisfies
\begin{equation}\label{NSomconv}
  \int_\Omega |\omega(x,t) - \alpha \omega^\chi(x,t)|\dd x \,=\, 
  \cO(t^{-\mu}\log t)\,, \qquad \hbox{as}\quad t \to \infty\,,
\end{equation}
where $\omega^\chi(x,t)$ is defined in \eqref{omegachi}, see
\cite{Ga2}.  In both convergence results \eqref{NSuconv},
\eqref{NSomconv}, one can replace the truncated Oseen vortex by the
original Lamb-Oseen vortex for which $\chi \equiv 1$, because the
additional error converges to zero like $\cO(t^{-1})$ as $t \to
\infty$.

\section{Axisymmetric vortex rings and filaments}\label{sec2}

When restricted to axisymmetric flows without swirl, the 
three-dimensional Navier-Stokes equations bear some similarity
with the two-dimensional situation considered in the previous
section. The only nonzero component of the vorticity vector 
satisfies Eq.~\eqref{Axiom}, which can be written in the equivalent 
form
\begin{equation}\label{Axiom2}
  \partial_t \omega_\theta + \partial_r(u_r\omega_\theta)  
  + \partial_z(u_z\omega_\theta) \,=\, \nu \Bigl((\partial_r^2 + \partial_z^2)
  \omega_\theta + \partial_r \frac{\omega_\theta}{r}\Bigr)\,.
\end{equation}
The analogy is most striking if one introduces the related 
quantity $\eta = \omega_\theta/r$, which satisfies the 
advection-diffusion equation
\begin{equation}\label{Axieta}
  \partial_t \eta + u\cdot\nabla\eta \,=\, \Delta \eta + 
  \frac{2}{r}\,\partial_r\eta\,,
\end{equation}
where $u\cdot\nabla = u_r\partial_r + u_z\partial_z$ and $\Delta =
\partial_r^2 + \frac{1}{r}\partial_r + \partial_z^2$. Equation 
\eqref{Axieta} is considered in the half-plane $\Omega = \{(r,z) \in 
\R^2\,|\, r > 0\,,~z \in \R\}$, with homogeneous Neumann boundary
conditions on $\partial\Omega$. 

It is clear from \eqref{Axieta} that $\eta(r,z,t)$ obeys the parabolic
maximum principle, and this provides a priori estimates on the
solutions which imply that the Cauchy problem for the axisymmetric
Navier-Stokes equations is globally well-posed, without any
restriction on the size of the initial data. The first results in this
direction were obtained by Ladyzhenskaya \cite{La} and by Ukhovskii
and Yudovich \cite{UY}, for finite energy solutions. Recently, it was
shown in \cite{GS1} that the vorticity equation \eqref{Axiom2} is
globally well-posed in the scale invariant space $L^1(\Omega, \D r\dd
z)$, equipped with the norm
\[
  \|\omega_\theta\|_{L^1(\Omega)} \,=\, \int_\Omega |\omega_\theta(r,z)|
  \dd r\dd z \,=\, \int_\Omega |\eta(r,z)|\,r\dd r\dd z\,.
\]
The proof follows remarkably the same lines as in the two-dimensional
case, and in particular uses very similar function spaces. Solutions
constructed in this way have infinite energy in general, but if the
initial vorticity $\omega_\theta^0$ decays somewhat faster at infinity
than what is necessary to be integrable, the velocity field becomes
square integrable for all positive times.

It is also possible to solve the Cauchy problem for Eq.~\eqref{Axiom2}
in the more general situation where the initial vorticity is a finite
measure. Global existence and uniqueness are established in \cite{GS1}
assuming that the total variation norm of the atomic part of the
initial vorticity is small compared to the viscosity parameter. The
general case is open to the present date, but interesting results have
been obtained for circular {\em vortex filaments}, which correspond to
the situation where the initial vorticity is a Dirac measure. The
resulting solutions can be considered as the analogue of the family of
Lamb-Oseen vortices in $\R^2$. These solutions cannot be written in
explicit form, but small-time asymptotic expansions can be computed
which involve the two-dimensional profiles $G$ and $v^G$ defined in
\eqref{Gdef}. If the circulation parameter $\Gamma$ is small compared
to the viscosity $\nu$, the results of \cite{GS1} imply the existence
of a unique global solution to \eqref{Axiom2} with initial vorticity
$\omega_\theta^0 = \Gamma \delta_{(\bar r, \bar z)}$, for any $(\bar
r, \bar z) \in \Omega$. For larger circulations, the following
existence result was recently established by Feng and \v{S}ver\'ak\:

\begin{prop}[\cite{FS}]\label{FSprop}
Fix $\Gamma > 0$, $(\bar r,\bar z) \in \Omega$, and $\nu > 0$. Then 
the axisymmetric vorticity equation \eqref{Axiom2} has a nonnegative 
global solution such that $\omega_\theta(t) \weakto \Gamma 
\delta_{(\bar r, \bar z)}$ as $t \to 0$. Moreover, this solution satisfies, 
for all $t > 0$, 
\begin{equation}\label{FSest}
  \int_\Omega \omega_\theta(r,z,t)\dd r\dd z \,\le\, \Gamma\,, \qquad
  \int_\Omega r^2 \omega_\theta(r,z,t)\dd r \dd z \,=\, 
  \Gamma\,\bar r^2\,.
\end{equation}
\end{prop}

Needless to say, the assumption $\Gamma > 0$ does not restrict the
generality, because the corresponding result for $\Gamma < 0$ can be
obtained by symmetry. Proposition~\ref{FSprop} is proved by a very
general approximation argument, which provides global existence without
any restriction on the size of the circulation parameter, but does not
imply uniqueness and does not give any precise information on the
qualitative behavior of the solution for short times. Using more
sophisticated techniques, a more accurate result can be established\:

\begin{thm}[\cite{GS2}]\label{GSmain}
Fix $\Gamma \in \R$, $(\bar r,\bar z) \in \Omega$, and $\nu > 0$. 
Then the axisymmetric vorticity equation \eqref{Axiom2} has a unique 
global mild solution $\omega_\theta \in C^0((0,\infty),L^1(\Omega) 
\cap L^\infty(\Omega))$ such that
\begin{equation}\label{GScond}
  \sup_{t > 0} \|\omega_\theta(t)\|_{L^1(\Omega)} \,<\, \infty~, 
  \qquad \hbox{and}\quad \omega_\theta(t) \weakto \Gamma 
  \,\delta_{(\bar r,\bar z)} \quad \hbox{as }t \to 0~.
\end{equation}
In addition, there exists a constant $C > 0$ such that the following
estimate holds\:
\begin{equation}\label{GSfinal}
  \int_\Omega \,\Bigl|\omega_\theta(r,z,t) - \frac{\Gamma}{4\pi\nu t}
  \,e^{-\frac{(r-\bar r)^2+(z-\bar z)^2}{4\nu t}}\Bigr|\dd r \dd z
  \,\le\, C \,|\Gamma|\,\frac{\sqrt{\nu t}}{\bar r}\,\log 
  \frac{\bar r}{\sqrt{\nu t}}~,
\end{equation}
as long as $\sqrt{\nu t} \le \bar r/2$.
\end{thm}

Since the existence of a global solution to \eqref{Axiom2} satisfying
\eqref{GScond} is already asserted by Proposition~\ref{FSprop}, the
main contributions of Theorem~\ref{GSmain} are the uniqueness of that
solution and its asymptotic behavior as $t \to 0$, as described in
\eqref{GSfinal}. The first step in the proof is a localization
estimate, which can be established using a Gaussian upper bound on the
fundamental solution of the ``linear'' equation \eqref{Axiom2} for
$\omega_\theta$, where the velocity field $u = (u_r,u_z)$ is
considered as given. It is found that, for any $\epsilon > 0$, there
exists a constant $C_\epsilon > 0$ such that
\begin{equation}\label{Gaussbound}
  |\omega_\theta(r,z,t)| \,\le\, \frac{C_\epsilon|\Gamma|}{\nu t}
  \,\exp\Bigl(-\frac{(r-\bar r)^2 + (z-\bar z)^2}{(4+\epsilon)
  \nu t}\Bigr)\,, \qquad (r,z) \in \Omega\,, \quad t > 0\,.
\end{equation}
Moreover $\int_\Omega \omega_\theta(r,z,t)\dd r\dd z$ converges to
$\Gamma$ as $t \to 0$. The second step consists in introducing
self-similar variables, in the spirit of \eqref{vwdef}. The 
rescaled vorticity $f$ and velocity $v$ are defined by
\[
  \omega_\theta(r,z,t) \,=\, \frac{\Gamma}{\nu t} \,f\Bigl(\frac{r 
  - \bar r}{\sqrt{\nu t}}\,, \frac{z -\bar z}{\sqrt{\nu t}}\,,\,t
  \Bigr)\,,\qquad u(r,z,t) \,=\, \frac{\Gamma}{\sqrt{\nu t}} 
  \,v\Bigl(\frac{r - \bar r}{\sqrt{\nu t}}\,, \frac{z - \bar z}{
  \sqrt{\nu t}}\,,\,t\Bigr)\,,
\]
and the following dimensionless quantities are introduced\:
\[
   R \,=\, \frac{r - \bar r}{\sqrt{\nu t}}~, \qquad  Z \,=\, 
  \frac{z - \bar z}{\sqrt{\nu t}}\,, \qquad 
  \epsilon \,=\, \frac{\sqrt{\nu t}}{\bar r}~, \qquad
  \alpha \,=\, \frac{\Gamma}{\nu}~.
\]
The evolution equation for the new function $f(R,Z,t)$ 
reads
\begin{equation}\label{fevol}
  t f_t + \alpha\Bigl(\partial_R(v_R f) + \partial_Z(v_Z f)\Bigr)
  \,=\, L f + \epsilon \partial_R\Bigl(\frac{f}{1+\epsilon R}\Bigr)\,,
\end{equation}
where as in \eqref{Ldef}
\[
  L \,=\, \partial_R^2 + \partial_Z^2 + \frac{R}{2}\,\partial_R + 
  \frac{Z}{2}\,\partial_Z + 1\,.
\]
Note that equation \eqref{fevol} lives in the time-dependent domain
where $1 + \epsilon R > 0$, but using the homogeneous Dirichlet
boundary condition one can extend the rescaled vorticity by zero
outside that domain and consider it as defined on the whole plane
$\R^2$. In the small time limit $\epsilon \to 0$, the system formally
reduces to the equation for perturbations around Oseen's vortex, which
was studied in detail in Section~\ref{ss1.1}, and the proof of
Theorem~\ref{GSmain} consists in showing that this intuition is indeed
correct. The Gaussian bound \eqref{Gaussbound} provides a uniform
control on the solution of \eqref{fevol} in the weighted space $X_t$
defined by the norm
\[
  \|f(t)\|_{X_t}^2 \,=\, \int_{1+\epsilon R > 0} f(R,Z,t)^2 \,e^{(R^2+Z^2)/4}
  \dd R \dd Z\,, \qquad t > 0\,,
\]
which coincides when $t = 0$ with the norm of the space $L^2(\infty)$
introduced in \eqref{def.L^2(m)}. A compactness argument, as in the
proof of Theorem~\ref{globconv}, can then be invoked to show that
$f(R,Z,t)$ necessarily converges to the Oseen vortex profile as 
$t \to 0$\:
\begin{equation}\label{fasym}
  \lim_{t \to 0} \|f(t) - \Gamma G\|_{X_t} \,=\, 0\,, \qquad 
  \hbox{where}\quad  G(R,Z) \,=\, \frac{1}{4\pi}\, e^{-(R^2+Z^2)/4}\,.
\end{equation}
The final step is an energy estimate which shows that, for some 
positive constants $C$ and $\delta$, 
\begin{equation}\label{fenergy}
  t\frac{\D}{\D t}\|f(t)\|_{X_t}^2 \,\le\, -\delta
  \|f(t)\|_{X_t}^2 + C\epsilon^2 |\log\epsilon|^2\,,
\end{equation}
when $t > 0$ is sufficiently small. The differential inequality
\eqref{fenergy} relies on the spectral properties of the linear
operator $L$ in the space $L^2(\infty)$, which were established in
Section~\ref{ss1.1}. As $X_t \hookrightarrow L^1(\R^2)$, it immediately
implies estimate \eqref{GSfinal} in Theorem~\ref{GSmain}. 
Moreover, a similar argument applied to the difference $f_1 - f_2$ 
of two solutions of \eqref{Axiom2} satisfying \eqref{GScond} 
leads to the conclusion that $f_1 \equiv f_2$, which yields 
uniqueness.

Theorem~\ref{GSmain} shows that the two-dimensional Lamb-Oseen vortex
naturally appears in the axisymmetric case too, where it describes the
short time behavior of solutions arising from vortex filaments as
initial data, see \eqref{GSfinal}. However, the long-time asymptotics
are very different in both situations, as can be seen from the
following result\:

\begin{prop}[\cite{GS1}]\label{GSmomprop}
Assume that the initial vorticity $\omega_\theta^0 \in L^1(\Omega)$ 
is nonnegative and has finite impulse\:
\begin{equation}\label{cM-def}
  \cI \,=\, \int_\Omega r^2 \omega_0(r,z)\dd r\dd z \,<\, \infty~.
\end{equation} 
Then the unique global solution of \eqref{Axiom2} satisfies
\begin{equation}\label{Mnon-V}
  \lim_{t \to \infty} \sup_{(r,z) \in \Omega} \,\left|t^2 \omega_\theta(r\sqrt{t},
  z\sqrt{t},t) - \frac{\cI}{16\sqrt{\pi}}\,r\,e^{-\frac{r^2+z^2}{4}}
  \right| \,=\, 0\,.
\end{equation}
In particular $\|\omega_\theta(t)\|_{L^\infty(\Omega)} = \cO(t^{-2})$
as $t \to \infty$. 
\end{prop}

Proposition~\ref{GSmomprop} applies in particular to the 
vortex rings constructed in Proposition~\ref{FSprop} and 
Theorem~\ref{GSmain}. It shows that the long-time asymptotics
are described, to leading order, by a self-similar solution 
of the {\em linearized} equation obtained by setting $u = 0$ 
in \eqref{Axiom2}. This is in sharp contrast with what happens
in the two-dimensional case. 

\section{Existence and stability of Burgers vortices}\label{sec3}

As was mentioned in the introduction, the Burgers vortex is a simple
but important model in fluid mechanics, describing the balance between
the dissipation due to the viscosity and the vorticity stretching
through the action of a background straining flow.  By rescaling
variables in a suitable manner (see, e.g. \cite{GW2}), one can 
assume without loss of generality that the rates of the linear 
strain in \eqref{defstrain} have the following form
\begin{equation}\label{defstrain.sec3}
 \gamma_1 \,=\,-\frac{1+\lambda}{2}\,,\qquad \gamma_2 \,=\, 
 -\frac{1-\lambda}{2}\,, \qquad \gamma_3 = 1\,.
\end{equation}
Here $\lambda \in [0,1)$ is a free parameter that represents the 
asymmetry of the strain, and the case $\lambda=0$ corresponds to 
an axisymmetric strain. The Burgers vortex with circulation $\alpha$ 
and asymmetry $\lambda$ is a two-dimensional stationary vorticity 
field of the form  $\oomega_{\lambda,\alpha}=(0,0,\omega_{\lambda,
\alpha})^\top$. In view of Eq.~\eqref{3Domstrain}--\eqref{LLdef}, this
means that the third component $\omega_{\lambda,\alpha}$ depends 
only on the horizontal variable $x = (x_1,x_2) \in \R^2$ and satisfies
the following elliptic problem in $\R^2$\:
\begin{equation}\label{eq.asBurgers}
  L_\lambda \omega - (v, \nabla )\omega \,=\, 0\,, \qquad v \,=\,
  K_{2D} * \omega\,, \qquad \int_{\R^2} \omega \dd x \,=\, \alpha\,,
\end{equation}
where $K_{2D}(x) = x^\perp/(2\pi|x|^2)$ and $L_\lambda$ is the two-dimensional 
differential operator defined by 
\begin{equation}\label{defLlambda}
  L_\lambda \,=\,\Delta + \frac{1+\lambda}{2}x_1\partial_1 
  + \frac{1-\lambda}{2}x_2\partial_2 + 1 \,=\,L + \lambda M\,,
  \qquad M \,=\,\frac{x_1}{2}\partial_1 - \frac{x_2}{2} \partial_2\,.
\end{equation}
When $\lambda=0$, Eq.~\eqref{eq.asBurgers} has the explicit solution
$\omega = \alpha G$, which is the classical axisymmetric Burgers
vortex \cite{Bu} with circulation $\alpha$. Note that $\alpha G$ is in
fact the unique solution of \eqref{eq.asBurgers} in the space
$L^1(\R^2)$, as can be deduced from Theorem~\ref{globconv}. Due to its
simple explicit expression, the axisymmetric Burgers vortex is often
used for comparison with experiments. However, the vortex tubes
observed in real flows or numerical simulations usually exhibit an
elliptical core region, rather than a circular one, because the local
strain due to the background flow is not axisymmetric in general. It
is therefore important to propose a model which takes into account the
asymmetry of the strain in an appropriate way, and allows one to
understand its influence on the shape of the vortex tubes. This
motivates the study of the Burgers vortex in the general case where
the asymmetry parameter $\lambda$ is nonzero \cite{MKO,PP2,RS}. In
that situation, solutions of \eqref{eq.asBurgers} cannot be written in
explicit form, and have to be constructed by a rigorous mathematical
argument. The aim of this section is to give an overview of the
mathematical results available by now about the existence of
asymmetric Burgers vortices (Section \ref{subsec3.1}) and their
stability with respect to two or three-dimensional perturbations
(Sections \ref{subsec3.2} and \ref{subsec3.3}).
 
\subsection{Existence and uniqueness of asymmetric Burgers 
vortices}\label{subsec3.1}

Since an explicit representation is no longer available for asymmetric
Burgers vortices, existence of such solutions is the first question to
address. One of the key observations is that, as the asymmetry
parameter $\lambda$ in \eqref{defstrain.sec3} is increased, the
localizing effect due to the linear strain becomes weaker in the $x_2$
direction. This phenomenon is illustrated by the shape of the
function
\begin{equation}\label{def.asG}
  \cG_\lambda (x) \,=\,\frac{\sqrt{1-\lambda^2}}{4\pi} 
  \,e^{-\frac{1+\lambda}{4}x_1^2 - \frac{1-\lambda}{4}x_2^2}\,,
  \qquad x \in \R^2\,,
\end{equation}
which solves the equation $L_\lambda \cG_\lambda \,=\,0$ in $\R^2$ with
$\int_{\R^2}\cG_\lambda\dd x = 1$. The form of $\cG_\lambda$ indicates
that asymmetric Burgers vortices, if they exist, still have a Gaussian
decay at infinity, but with a rate that becomes slower as $\lambda$
increases. Therefore, the function space $L^2(\infty)$ defined in
\eqref{def.L^2(m)} has to be modified in an appropriate way to allow
for a general asymmetry parameter $\lambda\in [0,1)$.  In view of
\eqref{def.asG} it is rather natural to introduce the function
\begin{align}\label{defG_lambda}
  G_\lambda (x) \,=\, \frac{1-\lambda}{4\pi}\,e^{-\frac{1-\lambda}{4}|x|^2}\,,
  \qquad x \in \R^2\,,
\end{align}
and the weighted $L^2$ spaces 
\begin{align}
  L^2 (\infty;\lambda) \,&=\,\Big \{f\in L^2 (\R^2)~\Big|~ 
  \|f\|_{L^2 (\infty;\lambda)}^2 := \int_{\R^2} | f(x)|^2 \frac{\dd x}{G_\lambda (x)} 
  <\infty  \Big \}\,,\\[2mm]
  L^2_0(\infty; \lambda) \,&=\,\Big \{ f\in L^2 (\infty; \lambda)~\Big|~ 
  \int_{\R^2} f\dd x= 0 \Big \}\,,
\end{align}
together with the associated weighted Sobolev space
$W^{1,2}(\infty;\lambda)$ equipped with the norm
\begin{align}
  \|f\|_{W^{1,2}(\infty; \lambda)} \,=\, \|\rho f\|_{L^2 (\infty;\lambda)} + 
  \|\nabla f\|_{L^2 (\infty; \lambda)}\,, \qquad \rho (x) = (1+|x|^2)^\frac12\,.
\end{align} 
When $\lambda=0$ these spaces are simply denoted by $L^2(\infty)$, 
$L^2_0(\infty)$, and $W^{1,2}(\infty)$, respectively.

In the analysis of the Burgers vortex, the linear operators $L_\lambda$ in
\eqref{defLlambda} and $\Lambda_f$ defined by
\begin{align}\label{defLambda_f}
  \Lambda_f  \,\omega \,=\, \big ( K_{2D}* f, \nabla \big ) \omega 
  + \big ( K_{2D} * \omega, \nabla ) f
\end{align}
play essential roles. The operator $\Lambda_f$ naturally appears as
the linearization of the quadratic term $(v,\nabla )\omega$ in
\eqref{eq.asBurgers} around a given vorticity profile $f$. Note that
$\Lambda_G$ is nothing but the operator $\Lambda$ defined in
\eqref{Lamdef}, but in the present section the general notation
$\Lambda_G$ will be preferred in order to emphasize the dependence
upon $G$. The operators $L_\lambda$ (when $\lambda = 0$) and
$\Lambda_f$ (when $f$ is radially symmetric) are invariant under
rotations about the origin in $\R^2$. It is thus natural to use polar
coordinates $(r,\theta)$ in the plane and to expand all functions in
Fourier series with respect to the angular variable $\theta$. In this
way, one can introduce the projections $\PP_n$ ($n \in \Z$) defined by
\begin{align}\label{defProjection}
  (\PP_n g)(r,\theta) \,=\,g_n (r) e^{i n \theta}\,, \qquad g_n (r) 
  \,=\,\frac{1}{2\pi} \int_0^{2\pi} g(r\cos\theta,r\sin\theta) 
  \,e^{-in\theta} \dd \theta\,, 
\end{align}
and the projected spaces $\PP_n X=\{\PP_n g\,|\,g\in X\}$ for any 
function space $X$ such as $L^2(\infty; \lambda)$. By construction
the projections $\PP_n$ commute with both operators $L$ and 
$\Lambda_G$. 

As long as the asymmetry parameter $\lambda$ lies in $[0,1)$, the
linear strain \eqref{defstrain.sec3} localizes the vorticity in the
horizontal directions, because $\gamma_1 < 0$ and $\gamma_2 < 0$.
Starting from this observation and relying on numerical calculations,
Robinson and Saffman \cite{RS} conjectured the existence of asymmetric
Burgers vortices, i.e., solutions to \eqref{eq.asBurgers}, for all
values of the parameters $\lambda \in [0,1)$ and $\alpha \in \R$, at
least in the regime where $\frac{\lambda}{1+|\alpha|}$ is small
enough. This fundamental question has been settled by now as follows.
 
\begin{thm}[Existence]\label{thm.subsec3.1.1} 
For all $\lambda\in [0,1)$ and all $\alpha\in \R$, there exists 
at least one asymmetric Burgers vortex $\omega_{\lambda,\alpha} 
\in L^2(\infty;\lambda)$ satisfying \eqref{eq.asBurgers}. 
\end{thm}

This existence result is established in \cite{GW3} when $0\le
\lambda \ll \frac12$ and $\alpha\in \R$, in \cite{GW2} when
$\lambda\in [0,1)$ and $|\alpha|\le \kappa (\lambda)\ll 1$, and in
\cite{M1} when $\lambda\in [0,\frac12)$ and $|\alpha|\ge R
(\lambda)\gg1$. Here $\kappa (\lambda)$ and $R(\lambda)$ are positive
numbers satisfying $\kappa(\lambda) \to 0$ as $\lambda \to 1$ and
$R(\lambda) \to \infty$ as $\lambda \to 1/2$. The proofs in 
\cite{GW2,GW3,M1} are based on a detailed analysis of some linearized 
operators, and existence of solutions to \eqref{eq.asBurgers} is 
established using the Banach fixed point theorem, which also 
provides uniqueness in a suitable class of functions, see
Theorem~\ref{thm.subsec3.1.3} below.  In contrast, the general
existence result for all $\lambda \in [0,1)$ and all $\alpha \in \R$ 
established in \cite{M2} relies on the Leray-Schauder fixed point 
theorem, which does not give any information about uniqueness.  

Since Burgers vortices are used to model tubular structures in
turbulent flows, it is highly interesting to study their
asymptotic shape as $|\alpha| \to \infty$ in the presence of
asymmetric linear strains. Numerical calculations by Prochazka and
Pullin \cite{PP1,PP2} and by Robinson and Saffman \cite{RS} indicate
that a large circulation $\alpha$ has a symmetrizing effect on the
vortex, so that the leading order of the flow in a bounded fluid
region is the axisymmetric Burgers vortex $\alpha G$, even when
$\lambda\ne 0$.  Using formal asymptotic expansions, Moffatt, Kida, and
Ohkitani \cite{MKO} explained this phenomenon analytically and
obtained for the Burgers vortex $\omega_{\lambda,\alpha}$ an
expression of the form
\begin{align}\label{eq.ansatz}
  \omega_{\lambda,\alpha} \,=\,\alpha G + W_\lambda + R_{\lambda,\alpha}\,, 
  \qquad \hbox{as}\quad |\alpha| \to \infty\,,
\end{align} 
where $W_\lambda$ is independent of $\alpha$ and $R_{\lambda,\alpha} 
= \cO(|\alpha|^{-1})$ as $|\alpha| \to \infty$. The second order 
term $W_\lambda$ in \eqref{eq.ansatz} is especially relevant, 
since it describes how the asymmetry of the background flow modifies 
the shape of the vortex. 

Inserting the formal expansion \eqref{eq.ansatz} into 
\eqref{eq.asBurgers} and using the cancellation $(K_{2D}*G, 
\nabla) G=0$, one obtains the relation
\begin{equation}\label{eq.prelim}
  L_\lambda \big ( \alpha G + W_\lambda + R_{\lambda,\alpha} ) \,=\,
  \alpha \Lambda_G W_\lambda + \cO(1)\,, \qquad \hbox{as}\quad 
  |\alpha| \to \infty\,.
\end{equation}
Since $L_\lambda G = (L+\lambda M) G = \lambda M G$, the terms 
proportional to $\alpha$ in both sides of \eqref{eq.prelim} are 
equal if and only if $W_\lambda = \lambda w_\infty$, where $w_\infty$ 
satisfies
\begin{align}\label{eq.w_infty}
  \Lambda_G\,w_\infty \,=\, M G \,. 
\end{align} 
Eq.~\eqref{eq.w_infty} is derived and studied numerically in \cite{MKO}, 
and discussed in \cite{GW3} within a rigorous functional framework. 
More precisely, it is shown in \cite[Proposition 3.1]{GW3} that 
there exists a unique $w_\infty$ in $\PP_2 W^{1,2}(\infty)
+ \PP_{-2} W^{1,2}(\infty)$ solving \eqref{eq.w_infty}, while
the expansion \eqref{eq.ansatz} has been mathematically verified 
as follows.

\begin{thm}[Large circulation asymptotics]\label{thm.subsec3.1.2}
For any $\lambda\in [0,1)$, there exists $R_0(\lambda)>0$ such that, 
for all $\alpha\in \R$ with $|\alpha|\ge R_0(\lambda)$, there exists 
an asymmetric Burgers vortex $\omega_{\lambda,\alpha} \in\PP^e 
L^2(\infty;\lambda)$ solving \eqref{eq.asBurgers} and 
satisfying
\begin{align}\label{est.thm.subsec3.1.2}
  \|\omega_{\lambda,\alpha} - \alpha G - \lambda w_\infty\|_{W^{1,2} (\infty;\lambda)} 
  \,\le\, \frac{\lambda C(\lambda)}{1+|\alpha|}\,.
\end{align}
Here $\PP^e$ is the even projection defined by $\PP^e = \oplus_{n\in \Z} 
\PP_{2n}$, and the constants $R_0(\lambda), C(\lambda)$ satisfy
\[
  \lim_{\lambda \to 1} R_0 (\lambda)  \,=\,\lim_{\lambda \to 1} 
  C (\lambda)  \,=\, \infty\,.
\]
\end{thm}

The conclusion of Theorem \ref{thm.subsec3.1.2} was first established
in \cite{GW3} assuming $0\le \lambda \ll \frac12$, in which case
estimate \eqref{est.thm.subsec3.1.2} holds in the stronger norm of
$W^{1,2} (\infty)$ and not just in $W^{1,2}(\infty;\lambda)$.  The
result was then extended to the larger range $\lambda\in [0,\frac12)$
in \cite{M1}, using again the space $W^{1,2}(\infty)$, and the general
case where $\lambda\in [0,1)$ was finally settled in \cite{M2}. The
basic strategy in \cite{GW3,M1,M2} is to construct a solution of
\eqref{eq.asBurgers} as a perturbation of the leading order
approximation $\alpha G + \lambda w_\infty$, in such a way that
estimate \eqref{est.thm.subsec3.1.2} holds. To explain this idea more
precisely, it is convenient to introduce the perturbation
$w^{(1)}=\omega_{\lambda, \alpha}-\alpha G-\lambda w_\infty$, which
has to solve the system
\begin{align}\label{eq.asBurgers.1}
\begin{split}
  \big ( L_\lambda -\alpha \Lambda_G - \lambda \Lambda_{w_\infty} \big ) 
  w^{(1)} \,&=\, (K_{2D}*w^{(1)}, \nabla )w^{(1)} + \lambda f_\lambda\,,  
  \qquad \int_{\R^2} w^{(1)} \dd x \,=\,0\,, \\
  f_\lambda \,&=\,-L_\lambda w_\infty + \lambda ( K_{2D}*w_\infty, \nabla) 
  w_\infty\,.
\end{split}
\end{align}
To show that $w^{(1)}$ is of order $\cO(|\alpha|^{-1})$ as $|\alpha| 
\to \infty$, the key observation is that the source term $f_\lambda$ in 
\eqref{eq.asBurgers.1} also belongs to the range of $\Lambda_G$. 
Indeed, a similar argument as in \cite[Proposition 3.1]{GW3} implies 
the existence of a unique $h_\lambda \in (I-P_0)\PP^e W^{1,2}(\infty)$ 
satisfying $\Lambda_G h_\lambda = f_\lambda$. Then $f_\lambda$ 
is decomposed as
\begin{align*}
  f_\lambda \,=\, \Lambda_G h_\lambda \,=\, -\frac{1}{\alpha} \big( 
  L_\lambda -\alpha \Lambda_G -\lambda \Lambda_{w_\infty} \big ) h_\lambda 
  + \frac{1}{\alpha} \big ( L_\lambda - \lambda \Lambda_{w_\infty} 
  \big ) h_\lambda\,,
\end{align*}
and Eq. \eqref{eq.asBurgers.1} can thus be reduced to the following
system for $w^{(2)}=w^{(1)} + \frac{\lambda}{\alpha} h_\lambda$\:
\begin{align}\label{eq.asBurgers.2}
\begin{split}
  \big ( L_\lambda -\alpha \Lambda_G - \lambda \Lambda_{w_\infty - 
  \frac{1}{\alpha} h_\lambda}  \big ) w^{(2)}  & \,=\, (K_{2D}*w^{(2)}, 
  \nabla )w^{(2)} + \frac{\lambda}{\alpha} F_{\lambda,\alpha} \,,  \qquad 
  \int_{\R^2}  w^{(2)} \dd x \,=\,0\,, \\
   F_{\lambda,\alpha} & \,=\,\big ( L_\lambda -\lambda \Lambda_{w_\infty} \big ) 
   h_\lambda + \frac{\lambda}{\alpha}  ( K_{2D} * h_\lambda, \nabla) h_\lambda \,.
\end{split}
\end{align}
It is clear that the source term $\frac{\lambda}{\alpha}
F_{\lambda,\alpha}$ is of order $\cO(\frac{\lambda}{|\alpha|})$ in
$L^2_0 (\infty;0)$ as $|\alpha|\to \infty$. Since
$\frac{\lambda}{\alpha}\Lambda_{h_\lambda}$ is a lower order
perturbation that becomes small in the regime where $|\alpha|\gg 1$,
the crucial step to establish \eqref{est.thm.subsec3.1.2} is to prove
the invertibility of the operator
\begin{align}\label{def.L.lambda.alpha}
  L_\lambda-\alpha \Lambda_G -\lambda \Lambda_{w_\infty}\,, 
  \qquad \hbox{in the space~} L^2_0(\infty; \lambda)\,,
\end{align}
together with a uniform estimate for its inverse when $|\alpha|\gg 1$. 

This is an easy task if $\lambda \in [0,\frac12)$ is small, because
one can then work in the space $L^2_0(\infty)$ instead of
$L^2_0(\infty; \lambda)$, and consider the operator in
\eqref{def.L.lambda.alpha} as a small perturbation of the more
familiar operator $L -\alpha \Lambda_G$, which has been thoroughly
studied in Section~\ref{sec1}. As was mentioned in \eqref{Lamasym},
the operator $\Lambda_G$ is skew-symmetric in $L^2 (\infty)$, namely
\[
  \langle \Lambda_G f, g \rangle_{L^2(\infty)} + \langle f, \Lambda_G g 
  \rangle_{L^2(\infty)} \,=\,0\,,\qquad \hbox{for all } f,g\in W^{1,2}(\infty)\,.
\]
Moreover, it follows from \eqref{Lconj} that $-L$ is a self-adjoint operator 
in $L^2(\infty)$ with compact resolvent, which satisfies the lower bound 
$-L \ge 0$ in $L^2 (\infty)$ and $-L \ge \frac12$ in $L^2_0(\infty)$. 
These two observations yield the uniform lower bound
\begin{align}
  \langle \big (-L + \alpha \Lambda_G \big )f, f\rangle_{L^2(\infty)} 
  \ge \frac12 \|f\|_{L^2 (\infty)}^2\,, \quad\qquad f\in L^2_0 (\infty)
  \cap D (L)\,,\label{est.L_Lambda}
\end{align}
which in turn implies that $\|(-L+\alpha\Lambda_G)^{-1}\|_{L^2_0
(\infty) \to L_0^2(\infty)} \le 2$ for all $\alpha \in \R$. A straightforward 
perturbation argument then gives a uniform estimate on $(-L_\lambda 
+ \alpha\Lambda_G + \lambda\Lambda_{w_\infty})^{-1}$ in $L_0^2(\infty)$ 
if $\lambda$ is sufficiently small, see \cite[Proposition 2.1]{GW3}. 

For larger values of $\lambda$, the term $\lambda\Lambda_{w_\infty}$
is not a small perturbation anymore, and the invertibility of the
operator $L_\lambda-\alpha\Lambda_G -\lambda \Lambda_{w_\infty}$ in
$L^2_0(\infty; \lambda)$ is not known for arbitrary $\alpha \in
\R$. However, if $|\alpha|$ is sufficiently large (depending on
$\lambda$), one can exploit the stabilizing effect that was already
discussed in Section~\ref{ss1.2} for the simpler operator $L-\alpha
\Lambda_G$, see Proposition~\ref{stabeff}.  For the modified operator
$L_\lambda - \alpha \Lambda_G$ with $\lambda\in [0,1)$, one has the
following estimates
\begin{align}\label{est.large.circulation}
\begin{split}
  & \lim_{|\alpha|\to \infty} \|(I-\PP_{0}) (-L_\lambda +\alpha 
  \Lambda_G)^{-1} \PP^e\|_{L^2_0(\infty;\lambda) \to 
  L_0^2(\infty;\lambda)} =0\,, \\ 
  & \lim_{|\alpha|\to \infty} \|\PP_{0} (-L_\lambda +\alpha 
  \Lambda_G)^{-1} (I-\PP_{0}) \PP^e\|_{L^2_0(\infty;\lambda) 
  \to L_0^2(\infty;\lambda)} =0\,, 
\end{split}
\end{align}
which are proved in \cite{M1} for $\lambda\in [0,\frac12)$ and in
\cite{M2} for arbitrary $\lambda\in [0,1)$. Roughly speaking, 
this means that the non-radially symmetric elements of $\PP^e 
L^2(\infty;\lambda)$ are strongly attenuated under the action 
of $(-L_\lambda +\alpha \Lambda_G)^{-1}$ when $|\alpha|$ is large. 
In addition, since the function $w_\infty$ belongs to $\PP_2
L^2(\infty) + \PP_{-2}L^2(\infty)$, one has the identity
\begin{equation}\label{idf}
  \PP_0 \Lambda_{w_\infty}\PP_0  f \,=\,0\, \qquad \hbox{hence}\quad
  \PP_0 \Lambda_{w_\infty} f \,=\, \PP_0 \Lambda_{w_\infty} (I-\PP_{0})f\,,
\end{equation}
for all $f\in \PP^e W^{1,2}(\infty;\lambda)$. Combining 
\eqref{est.large.circulation} and \eqref{idf}, it is possible 
to show that $\lambda \Lambda_{w_\infty}$ is a relatively 
small perturbation of $L_\lambda - \alpha\Lambda_G$ in $\PP^e 
L_0^2 (\infty; \lambda)$ if $|\alpha|$ is sufficiently large
(depending on $\lambda$), and this implies the invertibility 
of $L_\lambda - \alpha \Lambda_G +\lambda\Lambda_{w_\infty}$ in 
$\PP^e L^2_0 (\infty;0)$ and provides a uniform bound for the 
inverse when $|\alpha|\gg 1$.

When $\lambda\in [0,\frac12)$, the resolvent estimates
\eqref{est.large.circulation} also hold in the smaller space
$L^2_0(\infty)$, instead of $L_0^2(\infty,\lambda)$, and are
substantially easier to prove because one can then use the convenient
property that $\Lambda_G$ is skew-symmetric, see
\cite{GW3,M1}. However, when $\lambda \in [\frac12,1)$, the Burgers
vortex does not belong to $L^2(\infty)$ anymore, as can be seen from
the shape of the function $\cG_\lambda$ in \eqref{def.asG}. One is
then forced to work in a wider space such as $L^2(\infty;\lambda)$,
where the operator $\Lambda_G$ is no longer skew-symmetric.  The key
idea in \cite{M2} to overcome this difficulty is to construct
explicitly a bounded and invertible operator $T$ so that $\Lambda_G
T$, the right action of $T$ on $\Lambda_G$, becomes skew-symmetric.
With this skew-symmetrizer $T$, the equation $(L_\lambda -\alpha
\Lambda_G )w = f$ is written in the equivalent form
\[
  (L_\lambda -\alpha \Lambda_G T ) T^{-1} w \,=\,- L_\lambda (T-I )  
  T^{-1} w + f\,.
\]
Using the skew-symmetry of $\Lambda_G T$, one can show that the
operator $L_\lambda -\alpha \Lambda_G T$ satisfies resolvent bounds
similar to \eqref{est.large.circulation}, and the additional term
$L_\lambda(T-I)$ can be considered as a relatively small perturbation.
This implies the invertibility of $L_\lambda-\alpha \Lambda_G T +
L_\lambda (T-I)$, hence of $L_\lambda-\alpha\Lambda_G$. The argument
also yields uniform estimates on the inverse $(L_\lambda
-\alpha\Lambda_G)^{-1}$, which in turn make it possible to treat
$\lambda\Lambda_{w_\infty}$ as a relatively small perturbation when
$|\alpha|$ is sufficiently large (depending on $\lambda$), thus
concluding the proof of Theorem~\ref{thm.subsec3.1.2}.

The uniqueness of asymmetric Burgers vortices in a suitable class of
functions is also available for some range of parameters
$(\lambda,\alpha)$. 

\begin{thm}[Uniqueness]\label{thm.subsec3.1.3}
{\rm (i)} {\em Case $0\le \lambda\ll \frac12$ and $\alpha\in \R$\:} 
There exist $\lambda_0\in (0,\frac12)$ and $\tau_0>0$ such that, 
if $\lambda \in [0,\lambda_0]$ and $\alpha\in \R$, there exists at
most one asymmetric Burgers vortex $\omega_{\lambda,\alpha}$ in the
set $\big \{ f\in L^2 (\infty)~|~ \|f - \alpha G\|_{W^{1,2}(\infty)}
\le \tau_0 \big \}$. \\[1mm]
{\rm (ii)} {\em Case $0\le \lambda <1$ and $|\alpha|\ll 1$\:} For
all $\lambda\in [0,1)$ there exists $\kappa_0=\kappa_0(\lambda)>0$
such that, for any $\alpha\in \R$ with $|\alpha| \le \kappa_0$, there
exists at most one asymmetric Burgers vortex $\omega_{\lambda,\alpha}$
in the set $\big \{ f\in L^2 (\infty) ~|~ \|f-\alpha \cG_\lambda
\|_{L^2 (\infty)}\le \kappa_0 \big \}$. \\[1mm]
{\rm (iii)} {\em Case $0\le \lambda <1$ and $|\alpha|\gg 1$\:} For
all $\lambda\in [0,1)$ and all $\tau>0$ there exists
$R'(\lambda,\tau)\ge R_0 (\lambda)$ such that, for any $\alpha\in \R$
with $|\alpha|\ge R'(\lambda,\tau)$, there exists at most one
asymmetric Burgers vortex $\omega_{\lambda,\alpha}$ in the set $\big
\{ f\in \PP^e L^2(\infty,\lambda) ~|~ \|f- \alpha G - \lambda
w_\infty\|_{W^{1,2} (\infty;\lambda)} \le \tau \big \}$.  Here $R_0
(\lambda)$ is as in Theorem~\ref{thm.subsec3.1.2}, and
$R'(\lambda,\tau)$ satisfies $\displaystyle \lim_{\tau\to
\infty} R'(\lambda,\tau) \,=\,\infty$.
\end{thm}

The statement (i) of Theorem \ref{thm.subsec3.1.3} is proved in
\cite{GW3} using the uniform estimate \eqref{est.L_Lambda} for the
inverse of $L-\alpha\Lambda_G$ in $L^2_0(\infty)$, while (iii) is
established in \cite{M1,M2} using the stabilization effect at large
circulations described in \eqref{est.large.circulation}. The
uniqueness in the case (ii) is obtained in \cite{GW2} in the more
general framework of the polynomially weighted spaces $L^2(m)$. The
key point in the proof of (ii) is an estimate for the inverse
$L_\lambda^{-1}$ in $L^2_0 (m)$ when $m$ is large enough, which
enables to apply the Banach fixed point theorem when $|\alpha|$ is
sufficiently small. Remark that existence of asymmetric Burgers
vortices is also established in \cite{GW2,GW3,M1,M2} for all three
cases (i), (ii), and (iii) above, whereas uniqueness for the 
parameter regions not covered by Theorem \ref{thm.subsec3.1.3} is 
an interesting but difficult question, which is essentially open.
 
\subsection{Two-dimensional stability of asymmetric Burgers 
vortices}\label{subsec3.2}

In the parameter regions where existence and uniqueness have been
established, the next important issue is stability.  Since the Burgers
vortex itself is a two-dimensional vorticity field, it is possible to
study its stability within the class of purely two-dimensional flows,
and this is the point of view adopted in this subsection. The
axisymmetric case where $\lambda=0$ was already discussed in detail in
Section \ref{sec1}, hence the main focus here will be on the
asymmetric case $\lambda\ne 0$.

The evolution equations for the perturbations are obtained from system
\eqref{3Domstrain}, where $\nu=1$ and $\gamma_1, \gamma_2, \gamma_3$
are as in \eqref{defstrain.sec3}, by expanding the vorticity vector
$\oomega(\xx,t)$ around the stationary Burgers vortex
$\oomega_{\lambda,\alpha}(\xx) = (0,0,\omega_{\lambda,\alpha}(x))^\top$.  
When the vorticity field $\ww(\xx,t) = (0,0,w (x,t))^\top$ of the
perturbation is two-dimensional, the problem is reduced to the
following equations for the scalar function $w$\:
\begin{equation}\label{eq.pasBurgers}
\left\{\begin{array}{l}
  \partial_t w - \big (L_\lambda - \Lambda_{\omega_{\lambda,\alpha}} \big ) w  
  + (v, \nabla )w \,=\,0\,, \qquad v \,=\,K_{2D} * w\,, \qquad t>0 \,,
  \\[2mm]
   \displaystyle w|_{t=0} \,=\,w_0\,, \quad \int_{\R^2} w_0 \dd x \,=\,0\,.
\end{array}\right. 
\end{equation}
Here the operators $L_\lambda$ and $\Lambda_f$ are defined by
\eqref{defLlambda} and \eqref{defLambda_f}, respectively. As can be
expected, the properties of the linearized operator $L_\lambda -
\Lambda_{\omega_{\lambda,\alpha}}$ play a crucial role in the
stability analysis. It is not difficult to show that $L_\lambda$
generates a $C_0$-semigroup in the polynomially weighted space 
$L^2(m)$ for $m < \infty$, and an analytic semigroup in the 
Gaussian weighted space $L^2(\infty;\lambda)$. In fact, the 
semigroup $e^{t L_\lambda}$ has the following explicit representation
\[
  \big (e^{t L_\lambda} f \big) (x) \,=\, \frac{e^t}{4\pi \sqrt{a_{\lambda} 
  (t) a_{-\lambda} (t)}}\int_{\R^2}  \exp \bigg (-\frac{|x_1-y_1|^2}{4 
  a_{\lambda}(t)} - \frac{|x_2-y_2|^2}{4 a_{-\lambda}(t)}\bigg )\, 
  f (y_1 e^{\frac{1+\lambda}{2}t}, y_2 e^{\frac{1-\lambda}{2}t}) \dd y\,, 
\]
where $a_{\theta} (t) = (1 - e^{-(1+\theta)t})/(1+\theta)$. Since
$\Lambda_{\omega_{\lambda,\alpha}}$ is a relatively compact
perturbation of $L_\lambda$, the full linearized operator $L_\lambda -
\Lambda_{\omega_{\lambda,\alpha}}$ is also the generator of a $C_0$
(or analytic) semigroup, and the main concern is the long-time
behavior of that semigroup.  The following results have been
established in (essentially) the same parameter regions as in
Theorem~\ref{thm.subsec3.1.3}.

\begin{prop}[Linear stability]\label{thm.subsec3.2.1}
{\rm (i)} {\em Case $0\le \lambda\ll \frac12$ and $\alpha\in \R$\:} 
There exists $\lambda_1\in (0,\frac12)$ such that, for all
$\lambda\in [0, \lambda_1]$ and all $\alpha\in \R$, 
\begin{align}\label{est.thm.subsec3.2.1.1}
  \|e^{t ( L_\lambda - \Lambda_{\omega_{\lambda,\alpha}} )} f\|_{L^2 (\infty)} \le C 
  \|f\|_{L^2(\infty)} \,e^{-\frac{1-\lambda}{2} t}\,, \qquad t \ge 0\,,
\end{align}
for all $f \in L_0^2 (\infty)$. Here $C$ is a universal constant 
independent of $\lambda\in [0,\lambda_1]$ and $\alpha\in \R$.
\\[1mm]
{\rm (ii)} {\em Case $0\le \lambda <1$ and $|\alpha|\ll 1$\:}
For all $\lambda\in [0,1)$ there exists $\kappa_1 (\lambda)>0$ such 
that, if $|\alpha|\le \kappa_1(\lambda)$, then 
\begin{align}\label{est.thm.subsec3.2.1.2}
  \|e^{t ( L_\lambda - \Lambda_{\omega_{\lambda,\alpha}} )} f\|_{L^2 (m)} \le C 
  \|f\|_{L^2(m)} e^{-\frac{1-\lambda}{2} t}\,, \qquad t \ge 0\,,
\end{align}
for all $f \in L_0^2 (m)$, $m>3$. 
Here $C$ depends only on $\lambda\in [0,1)$ and $m$.\\[1mm]
{\rm (iii)} {\em Case $0\le \lambda<\frac12$ and $|\alpha|\gg 1$\:}
For all $\lambda\in [0,\frac12)$ there exists $R_1(\lambda)\ge R_0
(\lambda)$ such that, if $|\alpha|\ge R_1 (\lambda)$, then 
\begin{align}\label{est.thm.subsec3.2.1.3}
  \|e^{t ( L_\lambda - \Lambda_{\omega_{\lambda,\alpha}} )} f \|_{L^2(\infty)} \, \le 
  \, C\|f\|_{L^2 (\infty)} e^{- \frac{1-\lambda}{2}t}\,, \qquad t \ge 0\,,
\end{align}
for all $f \in L_0^2 (\infty)$.  
Here $C$ depends only on $\lambda$ and $\alpha$, while $R_0(\lambda)$ is
as in Theorem~\ref{thm.subsec3.1.2}.
\end{prop}

The statement (i) of Proposition~\ref{thm.subsec3.2.1} is proved in
\cite{GW3}, where $L_\lambda - \Lambda_{\omega_{\lambda,\alpha}}$ is
regarded as a small perturbation of the simpler operator $L-\alpha
\Lambda_G$ for which, as recalled in \eqref{est.L_Lambda}, the
stability estimate is obtained uniformly in $\alpha$ using the
skew-symmetry of $\Lambda_G$ in $L^2(\infty)$. The case (ii) follows
from the analysis developed in \cite{GW2}. In fact, as is mentioned in
the next subsection, the three-dimensional stability is the main
concern of \cite{GW2}, but the class of perturbations considered there
includes purely two-dimensional flows. In case (ii) the asymmetric
Burgers vortex $\omega_{\lambda,\alpha}$ is of order $\cO(|\alpha|)$ in
$L^2(\infty; \lambda)$, and the operator $L_\lambda -\Lambda_{\omega_{
\lambda,\alpha}}$ is handled as a small perturbation of
$L_\lambda$, for which complete information on the spectrum and the
associated semigroup is available. Case (iii) is treated in
\cite{M1}, using in an essential way the stabilizing effect at large
circulations described in \eqref{est.large.circulation}. The
restriction $\lambda\in [0,\frac12)$ in (iii) is due to the fact that,
when $\lambda \ge \frac12$, the operator $L_\lambda -
\Lambda_{\omega_{ \lambda,\alpha}}$ has to be analyzed in the space
$L^2(\infty,\lambda)$, where $\Lambda_G$ is no longer skew-symmetric.
So far this difficulty could not be overcome for the stability 
problem, although existence and uniqueness of Burgers vortices 
were established by constructing a suitable skew-symmetrizer, as 
explained in Section~\ref{subsec3.1}.

It should be emphasized here that, in all cases (i), (ii), and (iii),
the temporal decay estimate for the semigroup $e^{t(L_\lambda -
  \Lambda_{\omega_{\lambda,\alpha}})}$ involves the exponent
$-\frac{1-\lambda}{2}$, which is known to be optimal. Indeed, by
differentiating the identity $L_\lambda \omega_{\lambda,\alpha} -
(K_{2D}*\omega_{\lambda,\alpha},\nabla)\omega_{\lambda,\alpha} = 0$
with respect to $x_2$, one observes that $\partial_2
\omega_{\lambda,\alpha}$ is an eigenfunction of $L_\lambda -
\Lambda_{\omega_{\lambda,\alpha}}$ for the eigenvalue
$-\frac{1-\lambda}{2}$.  Numerical results due to Prochazka and Pullin
\cite{PP2} indicate that $-\frac{1-\lambda}{2}$ is actually the
largest eigenvalue of $L_\lambda - \Lambda_{\omega_{\lambda,\alpha}}$
in $L^2_0(\infty,\lambda)$ for any $\lambda\in [0,1)$, but a
mathematical proof of this conjecture is still missing, except 
in the three cases stated in Proposition~\ref{thm.subsec3.2.1}.

The semigroup $e^{t(L_\lambda - \Lambda_{\omega_{\lambda,\alpha}})}$ has 
standard parabolic smoothing properties. Nonlinear stability 
with respect to small initial perturbations can thus be obtained 
by analyzing the integral equation associated with 
\eqref{eq.pasBurgers}\:
\begin{align}
  w(t) \,=\,e^{t ( L_\lambda - \Lambda_{\omega_{\alpha,\lambda}})} w_0 - \int_0^t 
  e^{(t-s) (L_\lambda - \Lambda_{\omega_{\alpha,\lambda}})} \big ( K_{2D}* w(s), 
  \nabla \big) w(s) \dd s\,,
\end{align}
and applying the conclusions of Proposition~\ref{thm.subsec3.2.1}. 
This gives the following result\:

\begin{thm}[Local 2D stability]\label{thm.subsec3.2.2} {\rm (i)}
{\em Case $0\le \lambda \ll \frac12$ and $\alpha\in \R$\:} There exists
$\epsilon>0$ such that, for all $\lambda\in [0, \lambda_1]$ and all
$\alpha\in \R$, the following statement holds.  For all initial data
$w_0\in L_0^2(\infty)$ such that $\|w_0\|_{L^2(\infty)}\le \epsilon$, 
Eq.~\eqref{eq.pasBurgers} admits a unique solution $w\in
C^0([0,\infty); L_0^2 (\infty))$, which satisfies
\begin{align}
  \|w(t)\|_{L^2(\infty)} \, \le \, C \|w_0\|_{L^2 (\infty)} 
  \,e^{-\frac{1-\lambda}{2}t}\,,\qquad t \ge 0\,.
\end{align}
Here the constant $C$ is independent of $\lambda\in [0,\lambda_1]$ 
and $\alpha\in \R$, while $\lambda_1$ is as in 
Proposition~\ref{thm.subsec3.2.1}.
\\[1mm]
{\rm (ii)} {\em Case $0\le \lambda <1$ and $|\alpha| \ll 1$\:} For all
$\lambda \in [0,1)$ there exists $\epsilon =\epsilon (\lambda)>0$ such
that, for any $\alpha\in \R$ with $|\alpha|\le \kappa_1(\lambda)$, the
following statement holds.  For all initial data $w_0\in L_0^2 (m)$,
$m>3$, such that $\|w_0\|_{L^2(m)}\le \epsilon$, Eq.~\eqref{eq.pasBurgers} 
admits a unique solution $w \in C^0([0,\infty);L_0^2(m))$, which satisfies
\begin{align}
  \|w(t)\|_{L^2(m)} \, \le \, C \|w_0\|_{L^2 (m)} \,e^{- \frac{1-\lambda}{2}t}\,,
\qquad t\ge 0\,.
\end{align}
Here $C$ depends only on $\lambda$ and $m$, while $\kappa_1(\lambda)$ 
is as in Proposition~\ref{thm.subsec3.2.1}.
\\[1mm]
{\rm (iii)} {\em Case $0\le \lambda<\frac12$ and $|\alpha|\gg 1$\:} 
For all $\lambda\in [0,\frac12)$ and any $\alpha \in \R$ with 
$|\alpha|\ge R_1(\lambda)$, there exists $\epsilon=\epsilon(\lambda,
\alpha)>0$ such that the following statement holds. For all initial 
data $w_0 \in L_0^2 (\infty)$ with $\|w_0\|_{L^2(\infty)}\le \epsilon$,
Eq.~\eqref{eq.pasBurgers} admits a unique solution $w\in C^0([0,\infty);
L_0^2 (\infty))$, which satisfies
\begin{align}
  \|w(t)\|_{L^2(\infty)} \, \le \, C \|w_0\|_{L^2 (\infty)} 
  \,e^{- \frac{1-\lambda}{2}t}\,,\qquad t\ge 0\,.
\end{align}
Here $C$ depends only on $\lambda$ and $\alpha$, while $R_1(\lambda)$ 
is as in Proposition~\ref{thm.subsec3.2.1}.
\end{thm}

As in Proposition~\ref{thm.subsec3.2.1}, the statement (i) of Theorem
\ref{thm.subsec3.2.2} is proved in \cite{GW3}, while (ii) follows from
the results of \cite{GW2}. The case (iii) is obtained in \cite{M1}.
It should be emphasized here that the basin of attraction in the case
(i) is uniform in the circulation number $\alpha$, as a consequence 
of the linear estimate \eqref{est.thm.subsec3.2.1.1} in Proposition
\ref{thm.subsec3.2.1}.

\subsection{Three-dimensional stability of Burgers vortices}
\label{subsec3.3}

The stability analysis becomes more complicated when the perturbations
are three-dimensional, because the vorticity field is no longer a scalar
quantity, and vortex stretching terms already appear in the
linearized operator.  The problem is highly nontrivial even in the
axisymmetric case $\lambda=0$, where Rossi and Le Diz{\`e}s \cite{RL}
have shown that the linearized operator does not have any
eigenfunction with nontrivial dependence upon the vertical variable.
Numerical evidence of linear stability with exponential decay of the
perturbations was obtained by Schmid and Rossi \cite{SR}, but their
analysis also reveals the occurrence of short-time amplification for
generic solutions. A mathematical understanding of the underlying
mechanisms, leading to a rigorous explanation of these observations, 
is an important and challenging question, for which significant 
progress has been made in recent years.

Starting from the vorticity equation \eqref{3Domstrain}, with $\nu=1$
and $\gamma_1,\gamma_2,\gamma_3$ as in \eqref{defstrain.sec3}, it is
easy to write the evolution equation for perturbations $\ww=\oomega -
\oomega_{\lambda,\alpha}$, where $\oomega_{\lambda,\alpha}=(0,0,
\omega_{\lambda,\alpha})^\top$ is the Burgers vortex with circulation
$\alpha$. The result is
\begin{equation}\label{eq.p3DBurgers}
\left\{\begin{array}{l}
  \partial_t \ww - \big ( \cL_\lambda - 
  {\bf \Lambda}_{{\oomega}_{\lambda,\alpha}} \big ) \ww  
  + (\vv, \nabla ) \ww  - (\ww, \nabla ) \vv \,=\,0\,,
  \quad \vv \,=\,K_{3D} * \ww\,, \quad t>0 \,,\\[2mm]
  \displaystyle \ww|_{t=0} \,=\,\ww_0\,,
  \qquad \nabla\cdot \ww_0=0\,, \qquad \int_{\R^2} w_{0,3} 
  (x, x_3) \dd x \,=\,0\,, \quad x_3\in \R \,,
\end{array}\right.
\end{equation}
where $K_{3D}$ is the kernel of the Biot-Savart law \eqref{BS3D}.
Here the operator $\cL_\lambda$ is given by \eqref{LLdef} with 
$\nu=1$, namely 
\begin{align}\label{defLLambda.sec3}
  \cL_\lambda \,=\,
  \begin{pmatrix} 
  L_\lambda +\partial_3^2 -x_3\partial_3 - \frac{3+\lambda}{2} \\[1mm]
  L_\lambda +\partial_3^2 -x_3\partial_3 - \frac{3-\lambda}{2} \\[1mm]
  L_\lambda +\partial_3^2 -x_3\partial_3 
\end{pmatrix}\,,
\end{align}
and $L_\lambda$ is the two-dimensional differential operator 
\eqref{defLlambda}. On the other hand, the operator 
${\bf \Lambda}_{\oomega_{\lambda,\alpha}}$ is defined by
\begin{equation}\label{defLambda_oomega}
  \begin{split}
  {\bf \Lambda}_{\oomega_{\lambda,\alpha}} \ww \,=\, \null &(K_{3D}* 
  \oomega_{\lambda,\alpha}, \nabla ) \ww + (K_{3D}*\ww, \nabla ) 
  \oomega_{\lambda,\alpha} \\ &- (\ww, \nabla ) 
  K_{3D}* \oomega_{\lambda,\alpha} - (\oomega_{\lambda,\alpha},\nabla) 
  K_{3D}*\ww\,.
  \end{split}
\end{equation}
The divergence-free condition as well as the zero mass condition
$\int_{\R^2} w_3 (x, x_3) \dd x = 0$ are preserved under the 
evolution defined by \eqref{eq.p3DBurgers}. Note that, at least 
formally, a divergence-free vector field $\ww = (w_1,w_2,w_3)^\top$ 
always satisfies the identity $\frac{\dd}{\dd x_3} \int_{\R^2} w_3
(x, x_3)\dd x=0$. This means that the condition $\int_{\R^2} w_{0,3}(x,x_3) 
\dd x =0$ in \eqref{eq.p3DBurgers} is a natural requirement 
on the initial data, which does not restrict the generality; 
see \cite[Section 1]{GM1} for a detailed discussion. 

Since the Burgers vortex itself is essentially a two-dimensional flow,
it is natural to choose a functional setting that allows for purely
two-dimensional perturbations, and more generally for perturbations
which do not decay to zero as $|x_3| \to \infty$. For this purpose, 
the following function spaces are introduced in \cite{GW2,GM1}\:
\begin{align}\label{def.Xm}
  X(m)  \,=\,BC (\R; L^2(m) )\,, \qquad X_0(m) \,=\, BC(\R;L^2_0(m))\,, 
\end{align}
as well as $\X (m) \,=\,X(m) \times X(m)\times X_0(m)$. Here
$BC(\R;L^2(m))$ denotes the space of all bounded and continuous
functions from $\R$ into $L^2(m)$, which is a Banach space equipped
with the norm $\|\phi\|_{X(m)} = \sup_{x_3\in \R} \|\phi (\cdot,
x_3)\|_{L^2 (m)}$, and $BC(\R;L^2_0(m))$ is the closed subspace
defined in a similar way. Since the leading order term $\cL_\lambda$
in the evolution equation \eqref{eq.p3DBurgers} contains the dilation
operator $-x_3\partial_3$, one cannot expect that the solutions will
be continuous in time in the uniform topology of $\X(m)$. To restore
continuity in time, it is convenient to work in $X_{loc}(m)$, which is
the very same space $X(m)$ equipped with the weaker topology induced
by the countable family of seminorms $\|\phi\|_{X_n (m)} \,=\,
\sup_{|x_3|\le n} \|\phi (\cdot,x_3)\|_{L^2(m)}$, for $n\in \N$. For 
vector valued functions, the space $\X_{loc}(m)$ is defined 
in a similar way by endowing $\X(m)$ with the localized topology.  

Using these notations, the local stability results available 
so far can be summarized as follows. 

\begin{thm}[Local 3D stability]\label{thm.subsec3.3.2} 
{\rm (i)} For all $\lambda\in [0,1)$ and all $\mu\in
(0,\frac{1-\lambda}{2})$ there exist $\epsilon =\epsilon (\lambda)>0$
and $\kappa_2(\lambda,\mu)\in (0,\kappa_1 (\lambda)]$ such that, for
any $\alpha\in \R$ with $|\alpha|\le \kappa_2(\lambda,\mu)$, the
following statement holds.  For all initial data $\ww_0\in
\X(m)$, $m>3$, with $\nabla \cdot \ww_0=0$ and $\|\ww_0
\|_{\X(m)}\le \epsilon$, Eq. \eqref{eq.p3DBurgers} admits a
unique solution $\ww\in L^\infty (\R_+; \X(m))\cap
C^0([0,\infty); \X_{loc} (m))$, which satisfies
\begin{align}\label{est.thm.subsec3.3.2.1}
  \|\ww (t)\|_{\X(m)} \le C \|\ww _0\|_{\X(m)} \,e^{-\mu t}\,, 
  \qquad t\ge 0\,.
\end{align}
Here $C$ depends only on $\alpha$, and $\kappa_1 (\lambda)$ is 
as in Proposition~\ref{thm.subsec3.2.1}.\\[1mm]
{\rm (ii)} Let $\lambda = 0$. 
For all $m\in (2,\infty]$ and all $\alpha\in \R$ there exists
$\epsilon=\epsilon (m,\alpha)>0$ such that the following statement
holds.  For all initial data $\ww_0 \in \X(m)$ with $\nabla
\cdot \ww_0=0$ and $\|\ww_0\|_{\X(m)}\le \epsilon$,
Eq. \eqref{eq.p3DBurgers} admits a unique solution $\ww \in L^\infty
(\R_+; \X(m))\cap C^0([0,\infty); \X_{loc} (m))$, which satisfies
\begin{align}\label{est.thm.subsec3.3.2.2}
  \|\ww (t)\|_{\X(m)} \le C \|\ww_0\|_{\X(m)} \,e^{-\frac{t}{2}}
  \,,\qquad t\ge 0\,.
\end{align}
Here $C$ depends only on $m$ and $\alpha$.
\end{thm}

The statement (i) of Theorem \ref{thm.subsec3.3.2} is proved in
\cite{GW2}, using estimates on the semigroup $e^{t\cL_\lambda}$
generated by the operator $\cL_\lambda$. In view of 
\eqref{defLLambda.sec3}, one has the representation
\begin{align}\label{def.sgLLambda}
  e^{t\cL_\lambda} \ww \,=\,\bigg ( e^{-\frac{3+\lambda}{2} t} 
  e^{t S_\lambda} w_1\,,\,  e^{-\frac{3-\lambda}{2} t} e^{t S_\lambda} w_2\,,\,  
  e^{t S_\lambda} w_3 \bigg )^\top \,,
\end{align}
where $S_\lambda$ is the differential operator defined by $S_\lambda =
L_\lambda + \partial_3^2 - x_3\partial_3$.  Since the operators $L_\lambda$ 
and $\partial_3^2 - x_3\partial_3$ act on different variables, 
it is possible to obtain the following explicit formula
\begin{align}\label{def.sgSLambda}
  \big ( e^{t S_\lambda} f \big )  (\xx ) \,=\,\frac{1}{\sqrt{4\pi a_1 (t)}} 
  \int_\R  \exp \big ( - \frac{|x_3 e^{-t} - y_3|^2}{4a_1 (t)} \big ) 
  \,\bigg ( e^{t L_\lambda} f \,  (\cdot, y_3) \bigg ) (x) \dd y_3\,,
\end{align}
where $a_1(t) = (1-e^{-2t})/2$ and $e^{t L_\lambda}$ is the two-dimensional 
semigroup encountered in Section~\ref{subsec3.2}. Useful estimates for
the semigroup $e^{t S_\lambda}$ in $X(m)$ are established in
\cite{GW2}, together with elementary spectral properties of the
generator $S_\lambda$. Note that it is possible to take 
$\mu=\frac{1-\lambda}{2}$ in estimate \eqref{est.thm.subsec3.3.2.1},
as can be shown using some arguments borrowed from \cite{GM1}. 

The statement (ii) of Theorem \ref{thm.subsec3.3.2} is established in
\cite{GM1}.  Remarkably, as in the two-dimensional case, the local
stability holds for all values of the circulation number $\alpha$, and
moreover the rate of convergence $e^{-\frac{t}{2}}$ is uniform in
$\alpha$.  Although the result of (ii) is stated only in the purely
axisymmetric case $\lambda=0$, by a standard perturbation argument it
is also possible to prove local stability of the asymmetric Burgers
vortex $\oomega_{\lambda,\alpha}$ if the asymmetry parameter $\lambda$
is sufficiently small, depending on $|\alpha|$. The proof of (ii) in
\cite{GM1} is based on the analysis of the linearized operator
$\cL-\alpha {\bf\Lambda}_\GG$ and its associated semigroup, where
$\cL$ and $\alpha{\bf\Lambda}_\GG$ are shorthand notations for the
operators $\cL_\lambda$ and ${\bf\Lambda}_{\oomega_{\lambda,\alpha}}$,
respectively, when $\lambda=0$. Since ${\bf \Lambda}_\GG$ is a lower 
order perturbation it is not difficult to construct the semigroup 
$e^{t (\cL-\alpha{\bf\Lambda}_\GG)}$ in $\X(m)$, but the main 
problem is to control the long-time behavior. The following 
result is the key achievement of \cite{GM1}.

\begin{prop}[Axisymmetric linear stability]\label{thm.subsec3.3.1} 
For all $m\in (2,\infty]$ and all $\alpha\in \R$ one has
\begin{align}
  \|e^{t (\cL -\alpha {\bf \Lambda}_\GG)} \ff\|_{\X (m)}\le C e^{-\frac{t}{2}} 
  \|\ff\|_{\X(m)}\,,\qquad  t\ge 0\,,
\end{align}
for all ${\bf f}\in \X (m)$, where $C$ depends only on $m$ and
$\alpha$. Moreover, if $\nabla\cdot \ff =0$ then $\nabla\cdot e^{t
(\cL -\alpha {\bf \Lambda})} \ff =0$.
\end{prop}

The proof of Proposition~\ref{thm.subsec3.3.1} in \cite{GM1} is based
on two important observations\:

\begin{enumerate}[(I)]

\item As an effect of vortex stretching, the vertical derivatives 
of the velocity and vorticity fields decay exponentially as 
$t \to \infty$, so that the long-time asymptotics are governed
by a two-dimensional vectorial system. 

\item When restricted to two-dimensional solutions, the linearized
operator $\cL-\alpha {\bf \Lambda}_\GG$ has symmetry properties 
which imply uniform stability for all values of the circulation
parameter. 

\end{enumerate}

In the rest of this section both mechanisms are explained in 
some detail for the more general semigroup $e^{t(\cL_\lambda - 
{\bf \Lambda}_{\oomega_{\lambda,\alpha}})}$, where $0 \le \lambda < 1$. 
Proposition~\ref{thm.subsec3.3.1} is stated and proved in \cite{GM1}
in the axisymmetric case $\lambda = 0$, but the arguments are 
robust and can be used to establish linear stability in
the asymmetric case too.

Property (I) above is due to a very specific dependence of the 
operator $\cL_\lambda - {\bf \Lambda}_{\oomega_{\lambda,\alpha}}$ upon the 
vertical variable $x_3$. Indeed, using the definition in 
\eqref{defLLambda.sec3}, it is straightforward to verify that 
$[\partial_3, \cL_\lambda] = -\partial_3$, where $[A,B] = AB-BA$ denotes
the commutator of $A$ and $B$. Moreover, since the Burgers vortex
$\oomega_{\lambda,\alpha}$ is a two-dimensional stationary solution, 
one has $[\partial_3, {\bf \Lambda}_{\oomega_{\lambda,\alpha}}]=0$.  
At the level of the semigroup, these identities imply that
\begin{align}\label{est.kvertical}
  \partial_3^k\,e^{t(\cL_\lambda - {\bf \Lambda}_{\oomega_{\lambda,\alpha}})} 
  \,=\,e^{-kt} \,e^{t(\cL_\lambda - {\bf \Lambda}_{\oomega_{\lambda,\alpha}})} 
  \partial_3^k\,, \qquad t\ge 0\,, 
\end{align}
for all integer $k\in \N$. Since the semigroup $e^{t(\cL_\lambda - 
{\bf \Lambda}_{\oomega_{\lambda,\alpha}})}$ grows at most exponentially 
in time, at a rate that depends only on $\lambda$ and $\alpha$, 
Eq.~\eqref{est.kvertical} shows that the $k^{\rm th}$ order 
vertical derivative of any solution to the linearized equation
$\partial_t \ww =(\cL_\lambda - {\bf\Lambda}_{\oomega_{\lambda,\alpha}})
\ww$ decays exponentially as $t \to \infty$, if $k \in 
\N$ is large enough. By a simple interpolation argument, 
it follows that any expression involving at least one vertical 
derivative of the solution becomes negligible in the long-time
regime, which means that one can restrict the analysis to 
the two-dimensional vectorial system obtained by disregarding
the vertical dependence of all quantities under consideration. 
More precisely, in view of \eqref{defLambda_oomega}, the operator 
${\bf \Lambda}_{\oomega_{\lambda,\alpha}}$ can be decomposed as
\begin{align*}
  {\bf \Lambda}_{\oomega_{\lambda,\alpha}} \ww 
  \,=\,{\bf \Lambda}_{\oomega_{\lambda,\alpha}}^{(1)} \ww 
  + {\bf \Lambda}_{\oomega_{\lambda,\alpha}}^{(2)} \ww 
  - {\bf \Lambda}_{\oomega_{\lambda,\alpha}}^{(3)} \ww
  - {\bf \Lambda}_{\oomega_{\lambda,\alpha}}^{(4)} \ww\,,
\end{align*}
where, with the notations $\nabla_h=(\partial_1,\partial_2)^\top$
and $\ww_h=(w_1,w_2)^\top$,
\begin{align}
\begin{split}
  & {\bf \Lambda}_{\oomega_{\lambda,\alpha}}^{(1)} \ww \,=\,( K_{2D}* 
  \omega_{\lambda,\alpha} , \nabla_h ) \ww\,, \qquad\hspace{7pt}  
  {\bf \Lambda}_{\oomega_{\lambda,\alpha}}^{(2)} \ww \,=\,( K_{3D}* \ww , 
  \nabla ) \oomega_{\lambda,\alpha}\,,\\ 
  & {\bf \Lambda}_{\oomega_{\lambda,\alpha}}^{(3)} \ww \,=\,( \ww_h , 
  \nabla_h ) K_{3D}* \oomega_{\lambda,\alpha}\,, \qquad 
  {\bf \Lambda}_{\oomega_{\lambda,\alpha}}^{(4)} \ww \,=\,\omega_{\lambda,\alpha}  
  \partial_3 K_{3D}* \ww \,.
\end{split}
\end{align}
The discussion above motivates the following decomposition of 
${\bf \Lambda}_{\oomega_{\lambda,\alpha}}$ into 2D and 3D parts\:
\begin{align}\label{2D3Ddecomp}
\begin{split}
  {\bf \Lambda}_{\oomega_{\lambda,\alpha}} \ww & \,=\,
  \begin{pmatrix}
   (u_{\lambda,\alpha},\nabla_h )  \ww_h - (\ww_h, \nabla_h) u_{\lambda,\alpha} \\
   \Lambda_{\omega_{\lambda,\alpha}} w_3 
  \end{pmatrix} \, + \, {\bf R}_{\lambda,\alpha} \ww\,,\\
  {\bf R}_{\lambda,\alpha}  \ww & \,=\, \Big ( (K_{3D}* \ww )_h - K_{2D}*w_3 , 
  \nabla_h \Big) \oomega_{\lambda,\alpha}  \, - \,  
  {\bf \Lambda}_{\oomega_{\lambda,\alpha}}^{(4)} \ww \,.
\end{split}
\end{align}
To derive \eqref{2D3Ddecomp} one uses the fact that
$\oomega_{\lambda,\alpha} = (0,0,\omega_{\lambda,\alpha})^\top$ is a
two-dimensional vorticity field, so that
$K_{3D}*\oomega_{\lambda,\alpha} = (u_{\lambda,\alpha},0)^\top$, with
$u_{\lambda,\alpha} = K_{2D}*\omega_{\lambda,\alpha}$. The operator
$\Lambda_{\omega_{\lambda,\alpha}}$, defined in \eqref{defLambda_f}, 
is artificially produced by adding and subtracting the 
term $(K_{2D}*w_3,\nabla_h)\omega_{\lambda,\alpha}$ in the 
right-hand side. As is shown in \cite[Proposition 4.5]{GM1}, 
all terms in the operator ${\bf R}_{\lambda,\alpha}$ involve at least 
one derivative with respect to $x_3$, and hence play a negligible 
role in long-time asymptotics.  Therefore, the problem is now
reduced to the analysis of the simpler operator $\cL_{\lambda,\alpha}$
defined by
\begin{align*}
\begin{split}
  & \cL_{\lambda,\alpha} \ww   \,=\,\cL_\lambda \ww \,-\,  
  \begin{pmatrix} (u_{\lambda,\alpha},\nabla_h )  \ww_h - 
  (\ww_h, \nabla_h) u_{\lambda,\alpha}\\ \Lambda_{\omega_{\lambda,\alpha}} w_3 
  \end{pmatrix} 
  \, = \, \cA_{\lambda,\alpha} \ww \, + \,  (\partial_3^2 - 
  x_3 \partial_3 ) \ww\,, 
\end{split}
\end{align*} 
where 
\begin{align}\label{def.cL_lambda_alpha}
\begin{split}
  & \cA_{\lambda,\alpha} \ww \,=\,
  \begin{pmatrix} \cA_{\lambda,\alpha,h} \ww_h \\
  \cA_{\lambda,\alpha,3} w_3 \end{pmatrix} \,=\,
  \begin{pmatrix}
  \big ( L_\lambda - \frac{3+\lambda}{2}  \big ) w_1 - 
  (u_{\lambda,\alpha},\nabla_h ) w_1 + (\ww_h, \nabla_h) u_{\lambda,\alpha,1}\\[1mm]
  \big ( L_\lambda - \frac{3-\lambda}{2} \big ) w_2 -  (u_{\lambda,\alpha},
  \nabla_h )  w_2 + (\ww_h, \nabla_h) u_{\lambda,\alpha,2}\\[1mm]
  \big ( L_\lambda - \Lambda_{\omega_{\lambda,\alpha}} \big ) w_3 
  \end{pmatrix}\,.
\end{split}
\end{align}
The crucial observation here is that the (vectorial) operator
$\cA_{\lambda,\alpha}$ acts only on the horizontal variable,
so that the semigroup $e^{t\cL_{\lambda,\alpha}}$ generated by
$\cL_{\lambda,\alpha} = \cA_{\lambda,\alpha} + \partial_3^2 -
x_3\partial_3$ can be expressed in terms of the 2D semigroup
$e^{t\cA_{\lambda,\alpha}}$ in the same way as in
\eqref{def.sgSLambda}. As a consequence, the long-time behavior of the
semigroup $e^{t (\cL_\lambda - {\bf \Lambda}_{\oomega_{\lambda,\alpha}})}$ in 
$\X(m)$ can be deduced from the spectral analysis of the two-dimensional
operator $\cA_{\lambda,\alpha}$ acting on $\mathbb{L}^2(m) := 
L^2(m)^2\times L^2_0(m)$. This leads to the following 
criterion \cite{GM1}\:

\medskip\noindent
{\bf Stability criterion\:} {\em Let $\lambda\in [0,1)$, $\alpha\in \R$, 
and $m>2$. If the stability estimate 
\begin{align}\label{est.criterion}
  \|e^{t \cA_{\lambda,\alpha}}\|_{\mathbb{L}^2(m)\to \mathbb{L}^2(m)} \le 
  C e^{-\mu t}\,, \qquad  t\ge 0\,,
\end{align}
holds for some $\mu \in (0,\frac{1-\lambda}{2}]$, then $\|e^{t
(\cL_\lambda- {\bf \Lambda}_{\oomega_{\lambda,\alpha}})}
\|_{\X(m)\to \X(m)}\le C' e^{-\mu t}$ holds
for all $t\ge 0$.}

\bigskip The key observation (II) concerns the structure of the 2D
operator $\cA_{\lambda,\alpha}$. From the definition
\eqref{def.cL_lambda_alpha} it is apparent that the horizontal
component $\ww_h = (w_1,w_2)^\top$ and the vertical component $w_3$
are completely decoupled under the action of $\cA_{\lambda,\alpha}$.
Furthermore, the third component $\cA_{\lambda,\alpha,3} = L_\lambda
- \Lambda_{\omega_{\lambda,\alpha}}$ acting on $w_3$ is exactly the
linearized operator at the Burgers vortex considered in
Section~\ref{subsec3.2}, when only two-dimensional perturbations are
allowed. Proposition~\ref{thm.subsec3.2.1} (iii) thus provides 
the desired stability estimate for the semigroup generated 
by $\cA_{\lambda,\alpha,3}$, uniformly for all $\alpha \in \R$
if the asymmetry parameter $\lambda$ is small enough. 

One of the main contributions of \cite{GM1} is the analysis of the
horizontal component $\cA_{\lambda,\alpha,h}$, which also has
a nice structure that allows to obtain a stability estimate for all
$\alpha \in \R$, at least if $\lambda=0$. The argument is as
follows. Since the operator $\ww_h \mapsto (\ww_h,\nabla_h)
u_{\lambda,\alpha} -(u_{\lambda,\alpha},\nabla_h)\ww_h$ is a relatively compact
perturbation of the second order operator $L_\lambda$, a standard
perturbation argument reproduced in \cite[Proposition 3.4]{GM1} 
shows that the long-time behavior of the semigroup
$e^{t\cA_{\lambda,\alpha,h}}$ in $L^2(m)^2$ (with $m > 1$ 
sufficiently large) is determined by the eigenvalues of
the generator $\cA_{\lambda,\alpha,h}$ in a Gaussian weighted 
space such as $L^2(\infty;\lambda)^2$. As usual, when $\lambda\in 
[0,\frac12)$, one can use $L^2(\infty)^2$ instead of 
$L^2(\infty;\lambda)^2$. To locate the eigenvalues of
$\cA_{\lambda,\alpha,h}$, the following identities play 
a crucial role\:
\begin{align}\label{eq.identity.subsec.3h}
\begin{split}
  x_h \cdot \cA_{\lambda,\alpha,h} \ww_h \,=\,\null &(L_\lambda - 2) x_h 
  \cdot \ww_h-2 \nabla_h \cdot \ww_h -\lambda(x_1 w_1 - x_2 w_2) \\
  & - (u_{\lambda,\alpha}, \nabla_h )x_h \cdot \ww_h + 
  (\ww_h, \nabla_h ) x_h \cdot u_{\lambda,\alpha}\,,\\
  \nabla_h\cdot \cA_{\lambda,\alpha,h} \ww_h \,=\,\null &(L_\lambda -1 )
  \nabla_h\cdot \ww_h - (u_{\lambda,\alpha}, \nabla_h ) \nabla_h \cdot \ww_h \,.
\end{split}
\end{align}
Here the notation $x_h=x=(x_1,x_2)^\top$ is used. When $\lambda=0$
the two-dimensional velocity field $u_{\lambda,\alpha}  = \alpha v^G$
satisfies $x_h \cdot u_{\lambda,\alpha}=0$, hence the first identity in 
\eqref{eq.identity.subsec.3h} becomes substantially simpler. 
If $\ww_h \in L^2(\infty)^2 \cap D(L)$ is a nontrivial 
eigenfunction of $\cA_{\lambda,\alpha,h}$ with eigenvalue $\mu \in \C$, 
one has the obvious relations
\begin{align*}
  \mu \ww_h \,=\,\cA_{\lambda,\alpha,h}\ww_h\,, \qquad  \mu x_h \cdot 
  \ww_h \,=\,x_h \cdot \cA_{\lambda,\alpha,h} \ww_h\,,  \qquad  \mu 
  \nabla_h \cdot \ww_h \,=\,\nabla_h\cdot \cA_{\lambda,\alpha,h} \ww_h\,,
\end{align*}
which can be combined with \eqref{eq.identity.subsec.3h} to 
obtain valuable information on $\mu$. Indeed, assume for 
simplicity that $\lambda = 0$. If $\nabla_h \cdot \ww_h
\not\equiv 0$, the identity
\[
  \mu\nabla_h \cdot \ww_h \,=\, (L -1)\nabla_h\cdot \ww_h - 
  \alpha (v^G,\nabla_h ) \nabla_h \cdot \ww_h
\]
implies that $\Re(\mu) \le -3/2$, in view of the spectral 
properties of $L$ established in Proposition~\ref{Lspectrum}
and the fact that the operator $\omega \mapsto (v^G,\nabla)\omega$ 
is skew-symmetric in $L^2(\infty)$, see \eqref{Lamasym}.
If $\nabla_h \cdot \ww_h \equiv 0$ and $x_h \cdot \ww_h 
\not\equiv 0$, one has the relation
\[
  \mu x_h \cdot \ww_h \,=\, (L - 2)x_h \cdot \ww_h - 
  \alpha (v^G, \nabla_h)x_h \cdot \ww_h\,,
\]
which implies that $\Re(\mu) \le -2$ by the same argument. 
Finally, if $x_h \cdot \ww_h\equiv 0$, the eigenvalue 
equation reduces to 
\[
  \mu \ww_h \,=\, (L - {\textstyle\frac32})\ww_h - 
  \alpha (v^G, \nabla_h)\ww_h + \alpha \ww_h^\perp h\,, 
  \qquad h(x_h) \, = \, \frac{1}{2\pi |x_h|^2} (1- e^{-\frac{|x_h|^2}{4}})\,,
\]
and a simple energy estimate leads to the conclusion that $\Re(\mu)
\le -3/2$ in all cases. As a consequence, when $\lambda = 0$, one has
the desired stability estimate
\begin{equation}\label{last}
  \|e^{t\cA_{\lambda,\alpha,h}} \ww_h\|_{L^2 (m)^2} \,\le\, C_{m,\alpha} 
  \,e^{-\frac32 t} \|\ww_h\|_{L^2 (m)^2}\,, \qquad t \ge 0\,,
\end{equation}
for all $\alpha \in \R$, if $m > 1$ is sufficiently large. 
Note that the constant in \eqref{last} depends on $\alpha$ 
and becomes large as $|\alpha| \to \infty$, which may be 
related to the short time amplification phenomenon observed
numerically by Schmid and Rossi \cite{SR}. 
 
By a perturbation argument, it is easy to show that the stability
estimate also holds if the asymmetry parameter $\lambda$ is nonzero
and small, but it is unclear whether the smallness assumption on
$\lambda$ is uniform with respect to the circulation parameter
$\alpha$. This interesting question is answered affirmatively 
by Theorem \ref{thm.subsec3.2.2} if the perturbations are 
restricted to purely two-dimensional flows. In the general case 
what is missing so far is a precise information on the eigenvalues 
of the two-dimensional operator $\cA_{\lambda,\alpha,h}$, especially 
in the regime where $|\alpha|\gg 1$.

\section{Conclusion}\label{sec4}

As can be seen from the results reviewed in the previous sections, the
mathematical theory of viscous vortices has reached a certain level of
maturity, but many interesting questions remain open to the present
date. In the simple case of a single, axisymmetric, straight vortex
tube, there are explicit formulas for the vorticity and velocity
profiles, and the stability with respect to two-dimensional
perturbations is fully understood for all values of the total
circulation (Section~\ref{ss1.1}), including the large Reynolds number
limit where additional stabilization occurs (Section~\ref{ss1.2}).  In
presence of a non-axisymmetric strain, existence of stretched vortices
is known for all values of the circulation and asymmetry parameters
(Section~\ref{subsec3.1}), but uniqueness and stability results are
not completely satisfactory, except perhaps in the large circulation
limit where asymptotic symmetrization and stabilization are observed
(Section~\ref{subsec3.2}). When arbitrary three-dimensional
perturbations are allowed, local stability of the axisymmetric Burgers
vortex is well understood for all values of the total circulation
(Section~\ref{subsec3.3}), but less is known in the asymmetric case,
and the question is essentially open for the self-similar Lamb-Oseen
vortex, due to the lack of stretching in the vertical direction.

In real fluids, however, one usually observes the interaction of
several vortex tubes, none of which is perfectly straight, and
vorticity is also created near the boundaries. All these discrepancies
from the ideal situation considered above give rise to difficult
mathematical questions, which are essentially open. The rigorous
theory of curved vortex filaments in viscous flows is still in its
infancy, except perhaps in the axisymmetric case without swirl where
existence and uniqueness of vortex rings can be established
(Section~\ref{sec2}). Interaction of vortices has been studied so far
only in the weakly coupled regime where the distance between the
vortex centers is much larger than the size of the vortex cores
\cite{Ga0,Mar1,Mar2}. Stronger interactions, such as vortex merging (in
two dimensions) or reconnection of vortex tubes (in three dimensions),
play a crucial role in the dynamics of turbulent flows, but a rigorous
description of these phenomena seems completely out of reach. Finally,
there are no mathematical results yet concerning the interaction of
viscous vortices with rigid boundaries, although the existence of
self-similar vortices in two-dimensional exterior domains can be
established at least for small values of the circulation parameter
(Section~\ref{ss1.3}).

\section*{Cross references}

\begin{itemize}

\item Self-Similar Solutions to the Nonstationary Navier-Stokes Equations

\item Large Time Behavior of The Navier-Stokes Flow

\item Models and Special Solutions of the Navier-Stokes Equations

\item Inviscid Limit and Boundary Layer of the Navier-Stokes Flow

\end{itemize}






\end{document}